
\input amstex
\documentstyle{amsppt}
\input pictex
\magnification=\magstep1
\widestnumber\key{GHML}
\catcode`\@=11
\headline={\def\chapter#1{\chapterno@. }
\def\\{\unskip\space\ignorespaces}\headlinefont@
\rightheadline}
\def\logo@{}
\catcode`\@=\active
\NoBlackBoxes

\define\GlR#1{\text{\it Gl}(#1,\Bbb R)}
\define\glR#1{\frak {gl}(#1,\Bbb R)}
\define\GlC#1{\text{\it Gl}(#1,\Bbb C)}
\define\glC#1{\frak {gl}(#1,\Bbb C)}
\define\Gl#1{\text{\it Gl}(#1)}
\define\gl#1{\frak {gl}(#1)}
\define\Sl{\text{\it Sl}(2,\Bbb C)}
\define\SlF{\text{\it Sl}(2,\Bbb F)}
\define\slC{\frak {sl}(2,\Bbb C)}
\redefine\P#1{\Bbb P^{#1}}
\define\ov{\overline}
\define\und{\underline}
\define\om{\omega}
\define\th{\theta}
\define\Th{\Theta}
\define\w{\wedge}
\define\a{\text{\bf a}}
\define\J{\text{\bf J}}
\define\K{\text{\bf K}}
\define\bb{\text{\bf b}}

\define\pair#1#2#3{\left<#1,#2\right>_{#3}}
\define\ppair#1#2#3{\left<\!\left<#1,#2\right>\!\right>^{(#3)}}
\define\F{\frak F}
\redefine\O{\Cal O}
\define\C{\Cal C}
\define\N{\Cal N}
\define\U{\Cal U}
\define\V{\Cal V}
\define\W{\Cal W}
\define\Y{\Cal Y}
\define\Z{\Cal Z}
\define\g#1{{\frak g}_{#1}}

\define\Par{\par \bigpagebreak}
\define\[{$$}
\define\]{$$}
\define\tg#1{{\rm (#1)}}
\define\={-\!\!~}
\define\longto{\longrightarrow}
\define\ltop#1{\overset #1 \to \longrightarrow}
\define\math{{\tt MATHEMATICA }}

\topmatter
\title Exotic Holonomy on Moduli Spaces of Rational Curves
\endtitle

\author {Quo-Shin Chi} \footnotemark\\ \\ {Lorenz J Schwachh\"{o}fer}
\endauthor

\footnotetext{Supported in part by NSF grant DMS 9301060}
\address
\[ \matrix \format \l & \qquad \qquad \qquad \quad \r\\
\matrix \format \l\\
\text{Quo-Shin Chi}\\
\text{Department of Mathematics}\\
\text{Campus Box 1146}\\
\text{Washington University}\\
\text{St. Louis, Mo 63130, USA}\\
\text{\rm chi\@artsci.wustl.edu} \endmatrix &
\matrix \format \l\\
\text{Lorenz J Schwachh\"{o}fer}\\
\text{Max-Planck-Institut f\"{u}r Mathematik}\\
\text{Gottfried-Claren-Str. 26}\\
\text{53225 Bonn}\\
\text{Germany}\\
\text{\rm lorenz\@mpim-bonn.mpg.de}
\endmatrix \endmatrix \]
\endaddress


\abstract

Bryant \cite{Br} proved the existence of torsion free connections with
exotic holonomy, i.e. with holonomy that does not occur on the
classical list of Berger \cite{Ber}. These connections occur on moduli
spaces $\Y$ of rational contact curves in a contact threefold $\W$.
Therefore, they are naturally contained in the moduli space $\Z$ of all
rational curves in $\W$.

We construct a connection on $\Z$ whose restriction to $\Y$ is torsion
free. However, the connection on $\Z$ has torsion unless both $\Y$ and
$\Z$ are flat.

We also show the existence of a new exotic holonomy which is a certain
sixdimensional representation of $\Sl \times \Sl$. We show that every
regular $H_3$-connection (cf. \cite{Br}) is the restriction of a unique
connection with this holonomy.

\endabstract

\endtopmatter

\document

\hyphenation{ho-lo-no-my}
\hyphenation{ho-lo-no-mies}

\subheading{\S0 Introduction}

Since its introduction by {\` E}lie Cartan, the {\it holonomy} of a
connection has played an important role in differential geometry. Most
of the classical results are concerned with the holonomy of Levi Civita
connections of Riemannian metrics. In 1955, Berger \cite{Ber}
classified the possible irreducible Riemannian holonomies and much work
has been done since to study these holonomies and their applications.
See \cite{Bes} and \cite{Sa} for a historical survey and also \cite{J}
for more recent results.

At the same time, Berger also partially classified the possible
non-Riemannian holonomies of torsion free connections. However, his
classification omits a finite number of possibilities, which are
referred to as {\it exotic holonomies}. As of yet, the complete list of
exotic holonomies is still not known.

The incompleteness of Berger's list and therefore the existence of
exotic holonomies was shown by Bryant \cite{Br}. He investigated the
irreducible representations of $\SlF$, $\Bbb F = \Bbb R$ or $\Bbb C$.
For each $d \geq 1$, we can regard $\SlF$ as a subgroup $H_d \subseteq
\Gl{d + 1, \Bbb F}$ via the (unique) $(d + 1)$-dimensional irreducible
representation of $\SlF$ which will be described below. Moreover, if we
let $G_d \subseteq \Gl{d + 1, \Bbb F}$ be the centralizer of $H_d$,
then $G_d$ may be regarded as a representation of $\Gl{2, \Bbb F}$. For
$d \geq 3$, these representations do not occur on Berger's list of
possible holonomies and are therefore candidates for exotic
holonomies.

In his paper, Bryant showed that in the case $d = 3$ torsion free
connections with holonomies $H_3$ and $G_3$ do exist both if $\Bbb F =
\Bbb R$ and $\Bbb F = \Bbb C$. We shall refer to them as
$H_3$\=connections ($G_3$\=connections respectively). From now on, we
shall assume that $\Bbb F = \Bbb C$ unless stated otherwise.

$G_3$-structures occur naturally in the following way: let $\W$ be a
complex contact threefold and suppose there is a rational contact curve
$C$ in $\W$ such that the restriction of the contact line bundle $L|_C$
has degree $-3$. Then the moduli space $\Y$ of all close-by contact
curves carries a torsion free $G_3$-connection.

Conversely, every holomorphic torsion free $G_3$-connection is locally
equivalent to a connection on such a moduli space $\Y$ \cite{Br}. \Par

Before we proceed, let us briefly describe the irreducible
representations of $\Sl$, $\GlC2$ and $\GlC2 \times \Sl$.

For $n \in {\Bbb N}$, let $\V_n \subseteq {\Bbb C}[x,y]$ be the
$(n+1)$-dimensional subspace of homogeneous polynomials of degree $n$.
There is an $\Sl$-action on $\V_n$ induced by the transposed action of
$\Sl$ on ${\Bbb C}^2$, i.e. if $p \in \V_n$ and $A \in \Sl$ then \[
\text{$(A \cdot p)(x,y) := p(u,v)$ \quad with \quad $(u, v) = (x, y)
A$.} \]

Of course, this formula also describes an action of $\GlC2$ on $\V_n$.

The irreducible representations of $\GlC2 \times \Sl$ can be described
as follows: for $n, m \in {\Bbb N}$, we let $\V_{n,m} := \V_n \otimes
\V_m$. Then the action of $\GlC2 \times \Sl$ on $\V_{n,m}$ is defined
by \[ (A, B) \cdot (p \otimes q) := (A \cdot p) \otimes (B \cdot q) \]
with the actions of $\GlC2$ and $\Sl$ on $\V_n$ and $\V_m$ from above.
We define the subgroup $G_{n,m} \subseteq \Gl{\V_{n,m}}$ to be the
image of this representation. Also, we let $H_{n,m} \subseteq G_{n,m}$
be the image of $\Sl \times \Sl \subseteq \GlC2 \times \Sl$. In other
words, $H_{n,m} = G_{n,m} \cap \text{\it Sl} (\V_{n,m})$.

It is well known \cite{H} that these are complete lists of the
irreducible representations of $\Sl$, $\GlC2$ and $\GlC2 \times \Sl$
respectively.  \Par

Given a rational contact curve $C$ in $\W$ as above, it turns out that
its normal bundle $N_C \to C$ is equivalent to $\O(2) \oplus \O(2)$. By
Kodaira's Deformation Theorem \cite{K}, the moduli space $\Z$ of {\it
all} curves near $C$ is a smooth analytic manifold. Obviously, $\Y
\subseteq \Z$.

The tangent space $T_C\Z$ can be identified with $H^0(\O(2) \oplus
\O(2)) \cong \V_{1,2}$ in a natural way. Therefore, $\Z$ carries a
canonical $G_{1,2}$-structure.

The main objective of this paper is to investigate the correlation
between the $G_3$-structure on $\Y$ and the $G_{1,2}$-structure on
$\Z$. It had been conjectured in \cite{Br} that the latter structure is
torsion free. However, we show that almost the exact opposite is true.
Namely, we shall prove

\proclaim{Theorem 0.1} Let $\W$ be a complex contact threefold, let $C$
be a rational contact curve in $\W$ such that the restriction of the
contact line bundle $L|_C$ has degree $-3$, and let $\Z$ ($\Y$
respectively) be the moduli space of rational curves (rational contact
curves respectively) in $\W$ close to $C$. Then the canonical
$G_{1,2}$-structure on $\Z$ is torsion free if and only if the
$G_3$-connection on $\Y$ is flat.
\endproclaim

This means that we cannot in general expect the $G_{1,2}$-structure on
$\Z$ to be torsion free. However, we can make some statement about its
torsion.

\proclaim{Theorem 0.2} Let $\W$ be a complex threefold, let $C$ be a
rational curve in $\W$ whose normal bundle is equivalent to $\O(2)
\oplus \O(2)$ and let $\Z$ be the moduli space of curves in $\W$ close
to $C$, equipped with the canonical $G_{1,2}$-structure. Then there is
a subbundle $T \subseteq \Lambda^2 T^*\Z \otimes T\Z$ of rank four and
a unique $G_{1,2}$-connection on $\Z$ whose torsion is a section of
$T$.
\endproclaim

The point of Theorem 0.2. is that the rank of $\Lambda^2 T^*\Z \otimes
T\Z$ equals $90$, so $T$ has a large codimension. In other words,
Theorem 0.2. states that the $G_{1,2}$-structure on $\Z$ has ``very
little torsion''.

Theorem 0.1. raises the questions if there are any non-flat torsion
free $G_{1,2}$\=structures at all.

\proclaim{Theorem 0.3}
\roster
\item The holonomy of a torsion free $G_{1,2}$-connection is contained
in $H_{1,2}$. Thus, every torsion free $G_{1,2}$-structure admits a
(one-parameter family of) $H_{1,2}$-reductions.
\item A regular torsion free $H_{1,2}$-structure with full holonomy is
determined by three parameters. Thus, $H_{1,2}$-connections do exists,
and $H_{1,2}$ is therefore another exotic holonomy representation.
\endroster
\endproclaim

Comparing this result with Theorem 0.1. it follows that the torsion
free $G_{1,2}$-connections do not arise as moduli spaces of rational
curves in a contact threefold. However, we have the following
characterization of $H_3$-connections.

\proclaim{Theorem 0.4} Given a torsion free $G_{1,2}$-connection on a
sixfold $\Z$ and an imbedding $\Y \hookrightarrow \Z$ of a fourfold
$\Y$ such that the connection on $\Z$ restricts to a $G_3$-connection
on $\Y$, then the holonomy of this restriction is contained in $H_3$.

Conversely, if $\Y$ is a fourfold with a regular torsion free
$H_3$-connection, then there is a unique torsion free
$G_{1,2}$-connection on some sixfold $\Z$ and an imbedding $\Y
\hookrightarrow \Z$ such that the connection on $\Z$ restricts to the
connection on $\Y$.
\endproclaim

Regularity of an $H_3$-connection is a generic condition. For the exact
definition, see \cite{Br}.

As an interesting consequence, we conclude that for a given regular
$H_3$\=connection on $\Y$, there is more than one $G_{1,2}$-connection
extending the $H_3$-connection, but exactly one of these extensions is
torsion free. \smallskip

The calculations in this paper make extensive use of the representation
theory of $\Sl$ and $\Sl \times \Sl$, particularly an explicit version
of the Clebsch-Gordan formula. For details, we refer the reader to
\cite{H} and \cite{Br}.

The {\it Clebsch-Gordan formula} describes the irreducible
decomposition of a tensor product of irreducible $\Sl$-modules: \[ \V_m
\otimes \V_n = \bigoplus_{p = 0} ^{min(m,n)} \V_{m + n - 2 p}\]

A convenient tool to compute the decomposition of polynomials into
their irreducible components are the bilinear pairings \[ \matrix
\pair{\text{\ }}{\text{\ }}p : V_n \otimes V_m \longrightarrow
V_{n+m-2p}\\ \\ \matrix \displaystyle{\pair uvp = \frac{1}{p!}
\sum_{k=0}^p (-1)^k \binom pk \frac{\partial^p u}{\partial^k x
\partial^{p-k} y} \frac{\partial^p v}{\partial^{p-k} x \partial^k y}} &
\text{for} & u \in V_n, v \in V_m. \endmatrix \endmatrix \]

It can be shown that these pairings are $\Sl$-equivariant and therefore
are the projections onto the summands of the Clebsch-Gordan formula.

The Clebsch-Gordan formula for the irreducible representations of $\Sl
\times \Sl$ reads: \[ \V_{i_1,i_2} \otimes \V_{j_1,j_2} =
\bigoplus_{p_1 = 0}^{min(i_1,j_1)} \bigoplus_{p_2 = 0}^{min(i_2,j_2)}
\V_{i_1 + j_1 - 2 p_1,\ i_2 + j_2 - 2 p_2} \]

On these spaces, we define the pairings \[
\pair{\text{\ }}{\text{\ }}{p_1,p_2} : V_{i_1, i_2} \otimes V_{j_1,
j_2} \longrightarrow V_{i_1 + j_1 - 2 p_1,\ i_2 + j_2 - 2 p_2} \]
defined by \[ \pair {u_{i_1} \otimes v_{i_2}} {u_{j_1} \otimes v_{j_2}}
{p_1, p_2} := (\pair {u_{i_1}} {u_{j_1}} {p_1}) \otimes (\pair
{v_{i_2}} {v_{j_2}} {p_2}) \] with the pairing
$\pair{\text{\ }}{\text{\ }}{p_k}$ from above.

Again, it can be shown that these pairings are $\Sl \times
\Sl$-equivariant and therefore are the projections onto the summands of
the Clebsch-Gordan formula. \Par

In \S1, we define the notion of a $G$-structure and intrinsic torsion
which will be important in \S3. To demonstrate the usefulness of this
concept we include several examples. \smallskip

In \S2, we cite Kodaira's Deformation Theorem \cite{K} which states
that, under certain circumstances, the moduli space $\Z$ of compact
submanifolds of a given space $\W$ is itself a manifold. We construct a
natural $G$-structure on this moduli space where $G$ is the
automorphism group of the normal bundle of an element of $\Z$ in $\W$,
provided some stability condition (condition \tg A) is satisfied.

If $\Z$ is the moduli space of {\it rational curves} with positive
semistable normal bundle, then condition \tg A is satisfied and hence
we get a $G$-structure on $\Z$. In particular, if $\dim (\W) = 3$ we
obtain $G_{1,k}$-structures for some positive integer $k$.

We then construct a certain class of connections on this
$G_{1,k}$-structure, called {\it special connections}. These have the
property that the submanifolds $\Z_p \subseteq \Z$ with $p \in \W$,
consisting of all $C \in \Z$ which pass through $p$, are totally
geodesic. This yields some information about the torsion of special
connections. \smallskip

In \S3, we consider a {\it contact threefold} $\W$, and a rational
contact curve $C$. We let $\Y$ be the moduli space of rational contact
curves close to $C$ and $\Z$ be the moduli space of {\it all} curves
close to $C$. We then construct a $G_{k+1}$-structure on $\Y$ from the
$G_{1,k}$-structure on $\Z$, and show that every connection on $\Z$
restricts naturally to a connection on $\Y \subseteq \Z$.

If $k = 2$, i.e. if the normal bundle of each $C \in \Z$ is equivalent
to $\O(2) \oplus \O(2)$, then -- using special connections -- we show
that the intrinsic torsion of $\Z$ is a section of a certain rank four
bundle. This establishes Theorem 0.2. Also, the restriction of this
connection to $\Y$ is {\it torsion free}, hence we get a new proof that
the $G_3$-structure on $\Y$ is torsion free. \smallskip

In \S4, we set up the structure equations for a {\it torsion free}
$G_{1,2}$-connection. It turns out that the first Bianchi identity
forces the holonomy of such a connection to lie in $H_{1,2} \subseteq
G_{1,2}$, hence we consider the structure equations for torsion free
$H_{1,2}$-connections instead. These equations and their derivatives
are similar -- albeit more complex -- to the structure equations for
$H_3$-connections studied in \cite{Br}. In fact, methods similar to the
ones used in \cite{Br} allow us to solve the structure equations
explicitly. Their moduli space is then computed and we prove Theorem
0.3. \smallskip

Finally, in \S5 we put together the results from \S\S3 and 4. First of
all, we show that if $\Z$ is the moduli space of rational curves in a
threefold $\W$ and if the associated $G_{1,2}$-structure on $\Z$ is
torsion free, then $\Z$ must be locally symmetric. Second, we determine
those torsion free $G_{1,2}$-structures on $\Z$ which restrict to a
$G_3$-structure on some $\Y \subseteq \Z$. Since none of these
structures are locally symmetric, Theorem~0.1. follows. We also
demonstrate Theorem 0.4. using the classification of $H_3$-connections
from \cite{Br}.

We conclude by discussing some questions which our investigation
raises. Namely, we show that every $G_{1,2}$-structure whose torsion is
a section of the bundle $T \subseteq \Lambda^2 T^*\Z \otimes T\Z$ from
{\rm Theorem 0.2.} is locally equivalent to the moduli space of
rational curves in a fivefold $\Cal P$ which integrate a rank two
Pfaffian system on $\Cal P$.

For example, if $\Z$ is the moduli space of rational curves in a
threefold $\W$, we can achieve this by letting $\Cal P := \P{}T\W$ with
the canonical differential system \cite{EDS}, and identifying each
curve $C \subseteq \W$ with its canonical lift $\hat C \subseteq \Cal
P$.

Of course, there are many rank two Pfaffian systems which are not
locally equivalent to this contact structure on $\P{}T\W$ \cite{C}. An
interesting question is:

{\it Which rank two Pfaffian systems on a fivefold $\Cal P$ yield
torsion free $G_{1,2}$\=connections?}

Since the moduli space of torsion free $G_{1,2}$-connections is only
three dimensional by Theorem 0.3., those systems must be very special.
The answer to this question will also shed some light onto the
significance of $H_3$-connections. We shall pursue this analysis in a
sequel of the present paper. \smallskip

\subheading{\S1 $G$-structures and intrinsic torsion}

Let $M^n$ be a (real or complex) manifold of dimension $n$. Let $\pi:
\F \to M$ be the {\it coframe bundle} of $M$, i.e. each $u \in \F$ is a
linear isomorphism $u: T_{\pi(u)}M \tilde{\longrightarrow} \V$, where
$\V$ is a fixed $n$-dimensional (real or complex) vector space. Then
$\F$ is naturally a principal right $\Gl \V$-bundle over $M$, where the
right action $R_g: \F \to \F$ is defined by $R_g(u) = g^{-1} \circ u$.
The {\it tautological $1$-form} $\th$ on $\F$ with values in $\V$ is
defined by letting $\th(\xi) = u(\pi_*(\xi))$ for $\xi \in T_u\F$. For
$\th$, we have the $\Gl \V$-equivariance $R_g^*(\th) = g^{-1}\th$.

Let $G \subseteq \Gl \V$ be a closed Lie subgroup and let $\g{}
\subseteq \gl \V$ be the Lie algebra of $G$. A $G$-structure on $M$ is,
by definition, a $G$-subbundle $F \subseteq \F$. For any $G$-structure,
we will denote the restrictions of $\pi$ and $\th$ to $F$ by the same
letters. Given $A \in \g{}$ we define the vector field $A_*$ on $F$ by
\[ (A_*)_u = \frac{d}{dt} \left(R_{exp(tA)}(u)\right)|_{t=0}. \]

The vector fields $A_*$ are called the {\it fundamental vertical vector
fields} on $F$. It is evident that $\pi_*(A_*) = 0$ and thus $\th(A_*)
= 0$ for all $A \in \g{}$; in fact, $\{A_* | A \in \g{} \} =
\ker(\pi_*)$. Moreover, for $A, B \in \g{}$ it is well-known that
$[A_*, B_*] = [A,B]_*$.

Let $x \in M$ and $u \in \pi^{-1}(x)$. The Lie algebra $\g x := u^{-1}
\g{} u \subseteq \gl{T_xM}$ is independent of the choice of $u$, and
the union $\g F := \bigcup_x \g x$ is a vector subbundle of $T^*M
\otimes TM$.

Now we shall consider the {\it first Spencer map} $Sp: \V^* \otimes \gl
\V \to \Lambda^2\V^* \otimes \V$ which is defined by
skew-symmetrization of the first two factors of $\V^* \otimes \gl \V
\cong \V^* \otimes \V^* \otimes \V$. Since $\g{} \subseteq \gl \V$ we
may consider the restriction of $Sp$ to $\V^* \otimes \g{}$, and we
define $\g{}^{(1)}$ and $H^{0,2}(\g{})$ by requiring that the following
sequence be exact: \[ 0 \longto \g{}^{(1)} \longto \V^* \otimes \g{}
\ltop{Sp} \Lambda^2\V^* \otimes \V \ltop{pr} H^{0,2}(\g{}) \longto 0.
\tag1-1$$

In the same way, we can define vector bundles $\g F^{(1)}$ and
$H^{0,2}_F$ over $M$ via the exact sequence \[ 0 \longto \g F^{(1)}
\longto T^*M \otimes \g F \ltop{Sp} \Lambda^2T^*M \otimes TM \ltop{pr}
H^{0,2}_F \longto 0. \tag1-2$$

{}From now on, we will denote points in $M$ by $x$ and points in $F$ by
$u$. Moreover, $\xi, \xi'$ denote tangent vectors on $F$ and we let
$\und{\xi}_u = \pi_*(\xi_u)$, $\und{\xi}'_u = \pi_*(\xi'_u)$ etc.

A {\it connection} on $F$ is a $\g{}$-valued 1-form $\om$ on $F$
satisfying the conditions \[ \matrix \format \l & \l\\ \om(A_*) = A &
\text{\ \ for all $A \in \g{}$, and}\\ R_g^*(\om) = g^{-1} \om g &
\text{\ \ for all $g \in G$.} \endmatrix \tag1-3$$

Given a connection $\om$, its {\it torsion} $\Th$ is the $\V$-valued
2-form given by \[ \Th = d\th + \om \w \th. \tag1-4$$

{}From \tg{1-3} and \tg{1-4} it follows that there is a section $Tor$ of
$\Lambda^2 T^*M \otimes TM$ satisfying \[ \Th(\xi_u,\xi'_u) =
u\ (Tor(\und{\xi}_u,\und{\xi}'_u)) \text{\ \ \ for all $\xi_u, \xi'_u
\in T_uF$ and all $u \in F$.} \tag1-5$$

The connection $\om$ is called {\it torsion free} if $\Th = 0$.

Now let $\om'$ be another connection on $F$, and let $\Th'$ and $Tor'$
represent its torsion. From \tg{1-3} it follows that there is a section
$\alpha$ of the bundle $T^*M \otimes \g F$ such that \[ (\om' -
\om)(\xi_u) = u \ \alpha(\und{\xi}_u) \ u^{-1}. \tag1-6$$

{}From \tg{1-4} - \tg{1-6} we obtain for the torsion \[ (\Th' -
\Th)(\xi_u, \xi'_u) = u\ \left( \alpha (\und{\xi}_u) \cdot \und{\xi}'_u
- \alpha (\und{\xi}'_u) \cdot \und{\xi}_u \right),\] and hence, \[
\matrix \format \l & \l & \c\\ (Tor' - Tor)(\und{\xi}, \und{\xi}') & =
\alpha (\und{\xi}) \cdot \und{\xi}' - \alpha (\und{\xi}') \cdot
\und{\xi}\\ & = Sp(\alpha) (\und{\xi}, \und{\xi}') & \text{for all
$\und{\xi}, \und{\xi}' \in T_x M$}. \endmatrix \]

Thus, we conclude that \[ Tor' = Tor + Sp(\alpha). \tag1-7$$

This implies that the section $\tau := pr(Tor)$ of the bundle
$H^{0,2}_F$ is independent of the choice of $\om$ and therefore {\it
only depends on the $G$-structure $F$}.

\definition{Definition 1.1} Let $\pi: F \to M$ be a $G$-structure.
\roster
\item The vector bundle $H^{0,2}_F$ is called the {\it intrinsic
torsion bundle of $F$}.
\item The section $\tau$ of $H^{0,2}_F$ defined above is called the
{\it intrinsic torsion of $F$}.
\item $F$ is called {\it torsion free} if its intrinsic torsion $\tau = 0$.
\endroster
\enddefinition

The following Proposition is then immediate from \tg{1-7}.

\proclaim{Proposition 1.2} Let $\pi: F \to M$ be a $G$-structure and
let $\tau$ be its intrinsic torsion.
\roster
\item If $\sigma$ is any section of $\Lambda^2T^*M \otimes TM$ such
that $pr(\sigma) = \tau$ then there is a connection on $F$ whose
torsion section $Tor$ equals $\sigma$.
\item There is a torsion free connection on $F$ if and only if $F$ is
torsion free.
\item If $F$ is torsion free then there is a one-to-one correspondence
between torsion free connections on $F$ and sections of $\g F^{(1)}$.
In particular, if $\g{}^{(1)} = 0$ then the torsion free connection on
$F$ is unique.
\endroster
\endproclaim

We will give some examples for this concept.

\definition {Examples 1.3}
\roster
\item Let $G = O(p,q) \subseteq \Gl \V$ with $\V = \Bbb R^n$ and $n =
p+q$. A $G$-structure on $M^n$ is equivalent to a pseudo-Riemannian
metric on $M$ of signature $(p,q)$. One can show that $Sp: \V^* \otimes
{\frak o}(p,q) \to \Lambda^2 \V^* \otimes \V$ is an {\it isomorphism}.
Thus, $\g{}^{(1)} = 0$ and $H^{0,2}(\g{}) = 0$. Then Proposition 1.2.
implies that there is a {\it unique torsion free connection} on such a
$G$-structure. Of course, this reproves precisely the existence and
uniqueness of the Levi-Civita connection of a (pseudo-)Riemannian
metric. \cite{KN}
\item Suppose $n = 2m$ and let $G = \GlC m \subseteq \GlR n$. A
$G$-structure on $M^n$ is equivalent to an almost complex structure on
$M$. Then $H^{0,2} (\glC m) = \{ \phi \in \Lambda^2 (\Bbb C^n)^*
\otimes_{\Bbb R} \Bbb C^n\ |\ \phi(i x,y) = -i \phi(x,y)\}$. Moreover,
the intrinsic torsion is given by the {\it Nijenhuis tensor}. It is
well known that the vanishing of this tensor, i.e. the torsion freeness
of the $G$-structure, is equivalent to the integrability of the almost
complex structure. \cite{KN}
\item Suppose $n = 2m$ and let $G = \text{\it Sp}(m) \subseteq \GlR n$.
A $G$-structure on $M^n$ is equivalent to a non-degenerate 2-form $\om$
on $M$, i.e. $\om^m \neq 0$. One can show that $H^{0,2} (\frak{sp}(m))
= \Lambda^3 \Bbb R^n$ and that the intrinsic torsion is represented by
the 3-form $d\om$. Thus, the $G$-structure is torsion free if anf only
if $\om$ is a symplectic form.
\endroster
\enddefinition

{}From these examples it should become evident that for many naturally
arising $G$-structures the vanishing of the intrinsic torsion implies,
in some sense, the ``most natural integrability condition'' of the
underlying geometric structure.

\subheading{\S2 $G$-structures on moduli spaces of compact
submanifolds}

Let $\W$ be a complex manifold of (complex) dimension $d + r$.

\definition{Definition 2.1} By an {\it analytic family of compact
submanifolds of dimension $d$ of $\W$} we shall mean a pair $(\N, \Z)$
of a complex manifold $\Z$ and a complex analytic submanifold $\N$ of
$\W \times \Z$ of codimension $r$ with the property that for each $t
\in \Z$, the intersection $\N \cap (\W \times t)$ is a compact
connected submanifold of $\W \times t$ of dimension $d$.
\enddefinition

We call $\Z$ the {\it moduli space} of the family $(\N, \Z)$. The
canonical projections of $\N$ onto $\W$ and $\Z$ will be denoted by
$pr_1$ and $pr_2$ respectively.

\[ \beginpicture
\put {$\N$} [Bl] at 0 0
\put {$\W $} [Bl] at -20 -20
\put {${ }_{pr_1}$} [Bl] at -17 -2
\put {${ }_{pr_2}$} [Bl] at 15 -2
\put {$\Z$} [Bl] at  17 -20
\arrow <1mm> [0.55, 2] from -1 -1 to -12 -12
\arrow <1mm> [0.55, 2] from 7 -1 to 18 -12
\endpicture \]

For each point $t \in \Z$, we set \[ C_t \times t = \N \cap (\W \times
t). \]

We may identify $C_t \times t$ with $C_t$ and consider $\N$ as a {\it
family consisting of compact submanifolds $C_t$, $t \in \Z$, of $\W$}.

{}From now on, we shall use the notational convention $H^i(E):=
H^i(C,\O(E))$ for a vector bundle $E \to C$.

Let $N_t \to C_t$ be the normal bundle of $C_t$ in $\W$. There is a
natural imbedding $\eta_t: T_t\Z \hookrightarrow H^0(N_t)$ \cite{K} and
we shall use $\eta_t$ to regard $T_t\Z$ as a subspace of $H^0(N_t)$.

\definition{Definition 2.2} An analytic family $(\N, \Z)$ is called
{\it complete at $t \in \Z$} if $\eta_t$ is an isomorphism. $(\N, \Z)$
is called {\it complete} if it is complete at $t$ for all $t \in \Z$.
\enddefinition

We now state one of the most famous classical Theorems of the subject:

\proclaim{Kodaira's Deformation Theorem \cite{K}} Let $C \subseteq \W$
be a compact submanifold of $\W$ of dimension $d$. Let $N \to C$ be the
normal bundle of $C$ in $\W$. If $H^1(N) = 0$ then there exists a
complete analytic family $(\N, \Z)$ such that $C = C_{t_0}$ for some
$t_0 \in \Z$.
\endproclaim

Let $E \to C$ be a holomorphic vector bundle and denote the group of
equivalences of $E$ with itself by $Aut(E)$. Since each $\phi \in
Aut(E)$ induces an isomorphism $\hat{\phi}:~\O(E) \to \O(E)$, we obtain
a natural representation $\alpha_*: Aut(E) \to \Gl{H^*(E)}$.

\definition{Definition 2.3} An analytic family $(\N, \Z)$ is said to
{\it satisfy condition \tg A} if
\roster
\item"{\it (i)}" for any $t_1, t_2 \in \Z$ the normal bundles $N_{t_i}
\to C_{t_i}$, $i = 1,2$, are equivalent, and
\item"{\it (ii)}" the representation $\alpha_0: Aut(N_t) \to
\Gl{H^0(N_t)}$ is faithful and has closed image for all $t \in \Z$.
\endroster
\enddefinition

Consider now a complete analytic family $(\N, \Z)$ satisfying condition
\tg A. Let $E \to C$ be a vetor bundle which is equivalent to the
normal bundles $N_t \to C_t$ for all $t \in \Z$, and let $G:= Aut(E)$.
Let $\V := H^0(E)$, and let \[ \pi: \F \to \Z \] be the $\V$-coframe
bundle of $\Z$. Now consider the principal bundle $\pi: F \to \Z$ with
\[ F := \left\{ \left. \imath: \CD E @> >> N_t\\ @V VV @V VV \\ C @> >>
C_t \endCD \right| \matrix t \in \Z,\\ \imath \text{\ a bundle
equivalence}. \endmatrix \right\}. \tag2-1$$

We can define a bundle imbedding $\zeta: F \hookrightarrow \F$ (and
thereby justify the double use of the symbol $\pi$) as follows: given
$\imath \in F$, there is an induced isomorphism $\hat{\imath}: \V \to
H^0(N_t)$. Then $\eta_t \circ \hat{\imath}: \V \to T_t\Z$ is also a
linear isomorphism, thus $(\eta_t \circ \hat{\imath})^{-1}$ is a point
in $\F$. From condition \tg A it follows that the definition
$\zeta(\imath) := (\eta_t \circ \hat{\imath})^{-1}$ is indeed
one-to-one, and moreover, the image $\zeta(F)$ is a $G$-structure on
$\Z$. Identifying $F$ with $\zeta(F)$, we regard $F$ as a $G$-structure
on $\Z$, where the principal $G$-action on $F$ is given by $R_g(\imath)
:= \imath \circ g$.

If we set $\und{\theta} := \zeta^*(\theta)$, where $\theta$ is the
tautological form on $\F$, then \[ \und{\theta} (\xi_{\imath}) =
\zeta(\imath) ((\pi \circ \zeta)_*(\xi_{\imath})). \tag2-2$$ By a
slight abuse of language, we shall call $\und{\theta}$ the {\it
tautological form on $F$}. \Par

{\it Thus, for a complete analytic family $(\N, \Z)$ satisfying
condition \tg A, we have constructed an induced $G$-structure on the
moduli space $\Z$.} \Par

For the remainder of this section, $(\N,\Z)$ will stand for a complete
analytic family of {\it rational curves} satisfying condition \tg A,
i.e. we assume that $d=1$ and $C_t \cong \P1$ for all $t \in \Z$.

It is well-known that every $k$-dimensional vector bundle $E \to \P1$
satisfies $\O(E) \cong \O(m_1) \oplus \ldots \oplus \O(m_k)$ for some
integers $m_i$, $i = 1,\ldots,k$. Moreover, it is not hard to show that
$E$ satisfies condition {\it \tg{ii}} in Definition 2.3. if and only if
$m_i \geq 0$ for all $i$ and $m_1 + \ldots + m_k >0$. In this case, the
automorphism group decomposes as \[ \matrix G \cong \tilde{G} \times
\Sl & \text{with} & \tilde{G} \cong \GlC{n_1} \times \ldots \times
\GlC{n_l}, \endmatrix \] where the $n_i$'s are the multiplicities of
the $m_i$'s. \cite{GH}

An interesting question is to determine the intrinsic torsion of such a
$G$-structure or at least to understand its {\it vanishing}. To do
this, we will construct connections on $F$ and make some statements
about their torsion.

Let \[ \rho: \Sl \times \P1 \to \P1 \] denote the action of $\Sl$ on
$\P1$ by M{\"o}bius transformations. Let us fix once and for all the
reference point \[ x_0:= [0:1] \in \P1. \]

Consider $(\N,\Z)$ as before. Let \[ P := \left\{ \left. \und{\imath}:
\P1 \to C_t \right| t \in \Z,\ \und{\imath} \text{\ a biholomorphism}
\right\} \] be the {\it parameter space of $\Z$}. Then the obvious
projection \[ \pi_{P,\Z}: P \to \Z \] is a principal $\Sl$-bundle,
where the principal action is defined by $R_g(\und{\imath}) =
\und{\imath} \circ \rho(g)$. There is another projection \[ \pi_{F,P}:
F \to P \] which maps a bundle isomorphism $\imath: E \to N_t$ to the
underlying biholomorphism $\und{\imath}: \P1 \to C_t$. $\pi_{F,P}$
yields another principal bundle with structure group $\tilde{G} =
G/\Sl$. Finally, there is a projection \[ \pi_{P,\N}: P \to \N \] which
projects $\und{\imath}: \P1 \to C_t$ to $(\und{\imath}(x_0),t) \in \N$.
This projection yields a principal bundle whose structure group is the
stabilizer $G_{x_0} \subseteq G$.

Summarizing, we have the following commutative diagram:

\[ \beginpicture
\put {$F$} [Bl] at 0 80
\put {$P$} [Bl] at 0 50
\put {$\N$} [Bl] at 25 25
\put {$\Z$} [Bl] at 0 0
\put {$\W$} [Bl] at 50 0
\put {${ }_{pr_2}$} [Bl] at 7 22
\put {${ }_{pr_1}$} [Bl] at 40 20
\put {${ }_{\pi_{P,\N}}$} [Bl] at 13 45
\put {${ }_{\pi_{P,\Z}}$} [Bl] at -13 30
\put {${ }_{\pi_{F,P}}$} [Bl] at -13 70
\arrow <1mm> [0.55, 2] from 5 79 to 5 60
\arrow <1mm> [0.55, 2] from 5 49 to 5 10
\arrow <1mm> [0.55, 2] from 7 49 to 23 33
\arrow <1mm> [0.55, 2] from 23 24 to 9 10
\arrow <1mm> [0.55, 2] from 32 24 to 48 8
\endpicture \]

Now let us fix $t_0 \in \Z$ and let $C := C_{t_0}$. If we denote the
tangent and the normal bundle of $C$ by $\tau_C$ and $N_C$ respectively
then we have the exact sequence \[ 0 \to \tau_C \to T\W|_C \to N_C \to
0, \tag2-3$$ where $T\W$ denotes the holomorphic tangent bundle of
$\W$.

It is well known that this sequence splits. Also, from \tg{2-3} we get
the exact sequence \[ 0 \to Hom(N_C,\tau_C) \to Hom(T\W|_C,\tau_C) \to
Hom(\tau_C,\tau_C) \to 0, \tag2-4$$ which in turn induces the exact
sequence \[ 0 \longto H^0(Hom(N_C,\tau_C)) \longto
H^0(Hom(T\W|_C,\tau_C)) \ltop{\pi^*} H^0(Hom(\tau_C,\tau_C)) \tag2-5$$

Let \[ S_C := (\pi^*)^{-1} (id_{\tau_C}) \subseteq
H^0(Hom(T\W|_C,\tau_C)), \] where $id_{\tau_C}$ is regarded as an
element of $H^0( Hom( \tau_C,\tau_C))$. Even though $\pi^*$ need not be
surjective in general, $S_C$ is non-empty; namely, $S_C$ consists of
all splitting maps of the exact sequence \tg{2-3}. Therefore, $S_C$ is
an affine subspace of $H^0(Hom( T\W|_C,\tau_C))$ whose dimension equals
that of $H^0(Hom(N_C,\tau_C))$.

Condition \tg A implies that the exact sequences \tg{2-3} - \tg{2-5}
are independent of the choice of $t_0 \in \Z$, hence so is the
dimension of $S_C$. In fact, the union \[ S := \bigcup_{t \in \Z}
S_{C_t} \] forms an affine bundle over $\Z$, called the {\it
split-bundle of $\Z$}.

\proclaim{Lemma 2.4} Given a local section $\sigma: \Cal U \to S$, $t
\mapsto \sigma_t$ of the split-bundle $S$, let $P_{\Cal U} :=
\pi_{P,\Z}^{-1} (\Cal U)$. There is a unique holomorphic connection
$\hat{\sigma}$ on $\pi_{P,\Z}: P_{\Cal U} \to \Cal U$ such that \[
\rho_*(\hat{\sigma} (\xi), 0_{x_0}) = \und{\imath}_*^{-1}
(\sigma_t(\und{\xi})) \tag2-6$$ for all $\xi \in T_{\und{\imath}}P$,
where $t = \pi_{P,\Z} (\und {\imath})$ and $\und{\xi} = (pr_1 \circ
\pi_{P,\N})_* (\xi)$.
\endproclaim

\demo{Proof} First of all, note that equation \tg{2-6} is well defined:
$\hat{\sigma} (\xi) \in \slC$, and therefore both sides are contained
in $T_{x_0} \P1$.

Let $\und{\sigma}(\xi)$ be the right hand side of \tg{2-6}. Then
$\und{\sigma}$ is a holomorphic 1-form on $P_{\Cal U}$ with values in
$T_{x_0}\P1$. Moreover, it is easy to see that $\rho_*(A, 0_{x_0}) =
\und{\sigma} (A_*)$ for all $A \in \slC$, where $A_*$ denotes the
fundamental vector field corresponding to $A$.

We define a basis $\{A_1, A_2, A_3\}$ of $\slC$ by the equation \[
\pmatrix a & b\\ c & -a \endpmatrix = a A_1 + b A_2 + c A_3. \]

Clearly, $\rho_*(A_i, 0_{x_0}) = 0$ for $i = 1,2$. We define the
complex-valued 1-form $\hat{\sigma}_3$ by the equation $\hat{\sigma}_3
(\xi)\ \rho_*(A_3, 0_{x_0}) = \und{\sigma} (\xi)$, and let
$\hat{\sigma}_1 := \frak L_{A_2} (\hat{\sigma}_3)$ and $\hat{\sigma}_2
:= \frac 1 2 \frak L_{A_2} (\hat{\sigma}_1)$, where $\frak L$ denotes
the Lie derivative.

It is left to the reader to verify that the $\slC$-valued 1-form \[
\hat{\sigma} := \sum_{i} \hat{\sigma}_i A_i \] defines a connection
with the desired property and that this choice is unique. \qed
\enddemo

Geometrically, the interpretation of the connection $\hat \sigma$ is
the following. Suppose we have a local section $\sigma: \U \to S$ and a
curve $\gamma: I \to \U \subseteq \Z$ for some open set $I \subseteq
\Bbb C$, $0 \in I$. Then a horizontal lift $\ov \gamma: I \to P$ of
$\gamma$ can be regarded as a map $\Gamma: I \times \P1 \to \W$ such
that $\Gamma(t, \_)$ parametrizes $\gamma(t)$.

Given a parametrization $\und \imath: \P1 \to \gamma(0)$, we then
define $\Gamma$ uniquely by requiring that

\roster
\item $\Gamma(0,x) = \und \imath(x)$ for all $x \in \P1$, and
\item $\sigma_t (\frac \partial {\partial t} \Gamma(t,x)) = 0$ for all
$t \in I$ and all $x \in \P1$.
\endroster

It is then easy to verify that the $\Gamma$ thus determined is the
horizontal lift of $\gamma$ w.r.t. the connection $\hat \sigma$.

\definition{Definition 2.5} A holomorphic connection $\om$ on $\pi: F
\to \Z$ is called {\it special} if there exists a holomorphic section
$\sigma$ of the split-bundle $S \to \Z$, and a $\tilde{\g{}}$-valued
1-form $\tilde{\om}$ on $F$ such that \[ \om = \tilde{\om} +
\hat{\sigma} \] with $\hat{\sigma}$ as in Lemma 2.4. Here, we use the
decomposition $\g{} \cong \tilde{\g{}} \oplus \slC$.
\enddefinition

\proclaim{Proposition 2.6} Let $(\N,\Z)$ and $\pi: F \to \Z$ be as
before. Every $t_0 \in \Z$ has a neighborhood $\Cal U \subseteq \Z$
such that the restricted bundle $\pi: F_\U \to \U$ with $F_\U :=
\pi^{-1} (\U)$ admits a special connection.
\endproclaim

\demo{Proof} The proof is almost obvious: choose $\U$ sufficiently
small such that $\pi: F_\U \to \U$ admits a holomorphic connection
$\und{\om} = \tilde{\om} + \phi$ where $\tilde{\om}$ and $\phi$ are
holomorphic 1-forms with values in $\tilde{\g{}}$ and $\slC$
respectively. After shrinking $\U$ we may also assume that the
split-bundle $S$ admits a holomorphic section $\sigma$ over $\U$. Then
the form $\om := \tilde{\om} + \hat{\sigma}$ is a special connection.
\qed
\enddemo

An important characterization of special connections comes from the
following

\proclaim{Proposition 2.7} For every $p \in \W$, let $\Z_p := \{ t \in
\Z | p \in C_t \}$. If $\Z_p \neq \emptyset$ then $\Z_p$ is a smooth
submanifold of $\Z$ with $\text{\rm codim} (\Z_p) = \dim(\W) - 1$. The
tangent space of $\Z_p$ at $t \in \Z_p$ is $H^0 (C_t, \O(N_t) - p)
\subseteq H^0 (N_t) \cong T_t\Z$. Moreover, $\Z_p$ is totally geodesic
w.r.t. any special connection $\om$.
\endproclaim

\demo{Proof} The proof of the first two parts is left to the reader.

To show that $\Z_p$ is totally geodesic, let $t_0 \in \Z_p$ and pick a
biholomorphism $\und{\imath}_0: \P1 \to C_{t_0}$ such that
$\und{\imath}_0 (x_0) = p$. Then $\und{\imath}_0 \in \pi_{P,\Z}^{-1}
(t_0)$. We also pick a bundle isomorphism $\imath_0: E \to N_{t_0}$
such that $\imath_0 \in \pi_{F,P}^{-1} (\und{\imath}_0)$. Here, $E \to
\P1$ is a vector bundle which is isomorphic to $\O(N_{t}) \to C_{t}$
for all $t \in \Z$.

Let $\om = \tilde{\om} + \hat{\sigma}$ be a special connection on $F$
where $\sigma$ is a section of the split-bundle $S$. Let $I \subseteq
\Bbb C$ be an open neighborhood of $0$, and consider a geodesic
$\gamma: I \to \Z$ with $\gamma(0) = t_0$ and $\gamma'(0) \in
T_{t_0}\Z_p$. Let $\tilde{\und{\gamma}}: I \to P$ and $\tilde{\gamma}:
I \to F$ be the horizontal lifts of $\gamma$ to $P$ and $F$
respectively with $\tilde{\und{\gamma}}(0) = \und{\imath}_0$ and
$\tilde{\gamma}(0) = \imath_0$.

Define $\Gamma: I \times \P1 \to \W$ by $\Gamma(z,x) :=
\tilde{\und{\gamma}} (z) (x)$. Since $\tilde{\gamma}$ is horizontal and
thus, in particular, $\hat{\sigma} (\tilde{\und{\gamma}}'(z)) = 0$ for
all $z$, it follows that \[ \matrix \sigma_{\gamma(z)} \left( \frac
\partial {\partial z} \Gamma(z, x_0) \right) = 0 & \text{for all $z \in
I$.} \endmatrix \tag2-7$$

On the other hand, since $\gamma$ is a geodesic, we conclude from
\tg{2-2} that \linebreak $\zeta (\tilde{\gamma} (z)) (\gamma'(z)) \in
H^0 (E)$ is independent of $z$ and thus vanishes at $x_0$ for all $z$.
It follows that \[ \matrix \frac \partial {\partial z} \Gamma (z, x_0)
& \text{is {\it tangent to $C_{\gamma(z)}$} for all $z \in I$.}
\endmatrix \tag2-8$$

But \tg{2-7} and \tg{2-8} together imply that \[ \matrix \frac \partial
{\partial z} \Gamma (z, x_0) = 0 & \text{for all $z \in I$,} \endmatrix
\] and thus $\Gamma (z, x_0) = p$ for all $z \in I$. But this means
$\gamma(z) \in \Z_p$ for all $z$ and this completes the proof. \qed
\enddemo

This Proposition yields immediately

\proclaim{Corollary 2.8} If $\om$ is a special connection on $\Z$ and
$Tor: \Lambda^2 H^0 (N_t) \to H^0 (N_t)$ is its torsion then \[
Tor(\Lambda^2 H^0 (C_t,\O(N_t) - p)) \subseteq H^0 (C_t,\O(N_t) - p) \]
for all $p \in C_t$. \qed
\endproclaim

Since torsion is a {\it local} concept, Proposition 2.6. together with
Corollary 2.8 will allow us to make some assumptions about the
intrinsic torsion of $F$. This will be applied in the following
sections.

We shall also need one further property of special connections. Its
proof is immediate from \tg{2-2} and \tg{2-6}.

\proclaim{Proposition 2.9} Let $\und \V := H^0(\P1, \O(k) \oplus \O(k)
- x_0)$ be the space of global sections of $\O(k) \oplus \O(k)$ which
vanish at $x_0$ and let $\g{}' := \tilde {\g{}} \oplus \frak h_{x_0}
\subseteq \g{}$ where $\frak h_{x_0} \subseteq \slC$ is the
infinitesimal stabilizer of $x_0$ under $\rho$. Consider the projection
$pr_1 \circ \pi_{F,\N}: F \to \W$. If $\om$ is a special connection on
$F$ then \[ \ker(pr_1 \circ \pi_{F,\N})_* = \{ \xi \in TF\ |\ (\th +
\om) (\xi) \in \und \V \oplus \g{}' \}. \qed \]
\endproclaim

\subheading{\S3 Moduli spaces of rational contact curves}

In this entire section, we shall assume that $\dim(\W) = 3$ and that
$\W$ carries a holomorphic contact structure, i.e. a holomorphic line
bundle $L \subseteq T^*\W$ with the property that for every
non-vanishing local section $\kappa$ in $L$, the local 3-form $\kappa
\wedge d\kappa$ does not vanish anywhere.

By a standard notational ambiguity we will denote by $\O(n)$ both the
(unique) line bundle of degree $n$ over $\P1$ and the sheaf of germs of
holomorphic sections of this line bundle.

Let us first of all cite the following

\proclaim{Proposition 3.1} \cite{Br} Let $\W$ denote a complex contact
3-fold with contact line bundle $L \subseteq T^*\W$. Let $C \subseteq
\W$ be an imbedded rational contact curve, and suppose that $L|_C \cong
\O(-k-1)$ for some integer $k \geq 0$. Then
\roster
\item $N_C \cong \O(k) \oplus \O(k)$, where $N_C$ denotes the normal
bundle of $C$ in $\W$,
\item the moduli space $\Z$ of imbedded rational curves is smooth and
of complex dimension $2k + 2$ near $C$, and
\item the subspace $\Y \subseteq \Z$ of rational contact curves in $\W$
is a smooth analytic submanifold of $\Z$ of dimension $k + 2$.
\endroster
\endproclaim

For the remainder of this section we shall assume that $k>0$. It
follows that $\Z$ is a {\it complete analytic family of rational curves
satisfying condition} \tg A.

Let $E := \O(k) \oplus \O(k)$. Then \[ G := Aut(E) = \GlC2 \times \Sl,
\] where the first factor $\GlC2$ consists of those automorphisms which
fix the base space $\P1$, and the second factor $\Sl$ consists of
automorphisms which are induced by M{\"o}bius transformations of $\P1$.
As an $Aut(E)$-module, $T_C\Z \cong H^0(\O(k) \oplus \O(k)) \cong
\V_{1,k}$.\Par

Let $L^{\perp} \subseteq T\W$ be the 2-plane bundle anihilated by the
sections of $L$. For local sections $\xi, \xi'$ and $\kappa$ of
$L^{\perp}$ and $L$ respectively, the pairing $(\xi \w \xi', \kappa)
\mapsto d\kappa (\xi,\xi')$ is easily seen to be tensorial and
non-degenerate, hence induces a bundle isomorphism
$\Lambda^2(L^{\perp}) \ltop{\cong} L^*$. Also, we have the canonical
short exact sequence \[ 0 \longto L^{\perp} \longto T\W \longto L^*
\longto 0, \tag3-1$$ where $T\W$ is the holomorphic tangent bundle of
$\W$.

Now let $C \in \Y \subseteq \Z$. Since $C$ is a contact curve, we have
an inclusion $0 \to \tau \to L^{\perp}_{|C}$, where $\tau$ is the
tangent bundle of $C$, and from there it follows that \[ \tau \otimes
(L^{\perp}_{|C}/\tau) \cong \Lambda^2(L^{\perp}_{|C}) \cong L^*_{|C}
\cong \O(k+1).\]

Thus, since $\tau \cong \O(2)$, we must have $L^{\perp}_{|C}/\tau \cong
\O(k-1)$. From \tg{3-1} we also have the short exact sequence \[ 0
\longto L^{\perp}_{|C}/\tau \longto N_C \longto L^*_{|C} \longto 0,
\tag3-2$$ where $N_C$ denotes the normal bundle of $C$ in $\W$.

Since $H^1(L^{\perp}_{|C}/\tau) \cong H^1(\O(k-1)) = 0$, \tg{3-2}
induces the short exact sequence \[ 0 \longto H^0(L^{\perp}_{|C}/\tau)
\ltop{\imath} H^0(N_C) \ltop{pr} H^0(L^*_{|C}) \longto 0. \]

\proclaim{Lemma 3.2} Let $C \in \Y$ and $\xi \in T_C\Y \subseteq T_C\Z
\cong H^0(N_C)$. If $pr(\xi) \in H^0(L^*_{|C})$ vanishes at $p \in C$
of order at least two, then $\xi$ -- regarded as a section of $N_C$ --
vanishes at $p$.
\endproclaim

\demo{Proof}
Given $\xi \in T_C\Y$ as above, we pick a holomorphic curve $\gamma: I
\to \Y$ with $\gamma(0) = C$ and $\gamma'(0) = \xi$ where $I \subseteq
\Bbb C$ is an open neighborhood of $0$. Let $\Gamma: I \times \P1 \to
\W$ be a holomorphic map such that $\Gamma(t,\_)$ is a parametrization
of $\gamma(t)$ for all $t$. We may assume $\Gamma(0,x_0) = p$ with $x_0
:= [0:1] \in \P1$.

First, suppose that $\Gamma$ is a local biholomorphism from a
neighborhood of $(0, x_0)$ to $U \subseteq \W$. Then the holomorphic
vector fields \[ \matrix X:= \displaystyle{\frac{\partial
\Gamma}{\partial t}} (t,x) & \text{and} & Y:=
\displaystyle{\frac{\partial \Gamma}{\partial x}} (t,x) \endmatrix \]
are well defined on $U$.

Let $\kappa$ be a local contact form on $U$. Then we have \[
d\kappa(X_p,Y_p) = X_p(\kappa(Y)) - Y_p(\kappa(X)) - \kappa([X,Y]_p).
\] But all three terms on the right hand side vanish: the first one
vanishes because $Y$ is tangent to the contact curves $\gamma(t)$, thus
$\kappa(Y) \equiv 0$. The second one vanishes because -- by hypothesis
-- the function $\kappa(X): (U \cap C) \to \Bbb C$ vanishes of order
two at $p$. Finally, $[X,Y] = 0$ from the definition of $X$ and $Y$,
thus the third term vanishes as well.

The vanishing of $pr(\xi)$ at $p$ implies that $X_p \in L^{\perp}_p$,
hence $X_p, Y_p$ span $L^{\perp}_p$. But this together with
$d\kappa(X_p,Y_p) = 0$ implies that $(\kappa \w d\kappa)_p = 0$ which
is impossible.

Therefore, $\Gamma$ is {\it not} a local biholomorphism at $(0, x_0)$,
i.e. $\frac{\partial \Gamma}{\partial t} (0,x_0)$ must be {\it tangent
to} $C$. But this implies exactly that $\xi$ vanishes at $p$. \qed
\enddemo

\proclaim{Corollary 3.3} For every $C \in \Y$, the restriction $pr:
T_C\Y \to H^0(L^*_{|C})$ is an isomorphism.
\endproclaim

\demo{Proof} From Proposition 3.1. we know that $\dim(T_C\Y) = k+2 =
\dim(H^0(L^*_{|C}))$, and from Lemma 3.2. it follows that $\ker(pr)
\cap T_C\Y = 0$. \qed
\enddemo

Recall that the sequence \tg{3-2} is equivalent to \[ 0 \to \O(k-1) \to
\O(k) \oplus \O(k) \to \O(k+1) \to 0. \] It is easy to show that -- up
to equivalence -- the maps in this exact sequence are uniquely
determined. More specifically, one can show that for every $C \in \Y$,
there are bundle isomorphisms $\phi_C, \phi'_C$ and $\phi''_C$ such
that the diagram
\[ \CD 0 @>>> \O(k-1) @>>> \O(k) \oplus \O(k) @>>> \O(k+1) @>>> 0\\ @.
@V\phi'_CVV @V\phi_CVV @V\phi''_CVV\\ 0 @>>> L^{\perp}_{|C}/\tau @>>>
N_C @>>> L^*_{|C} @>>> 0 \endCD \tag3-3$$ commutes, and for the induced
commutative diagram
\[ \CD 0 @>>> \V_{k-1} @>\und{\imath}>> \V_{1,k} @>\und{pr}>> \V_{k+1}
@>>> 0\\ @. @V(\phi'_C)^*VV @V\phi_C^*VV @V(\phi''_C)^*VV\\ 0 @>>> H^0(
L^{\perp}_{|C}/ \tau) @>\imath>> H^0(N_C) @>pr>> H^0(L^*_{|C}) @>>> 0
\endCD,\tag3-4$$
we have $\und{\imath} (u_{k-1}) = x \otimes (y \cdot u_{k-1}) - y
\otimes (x \cdot u_{k-1})$, and $\und{pr}(u_1 \otimes v_k) = u_1 \cdot
v_k$ for all $u_i, v_i \in \V_i$. Here, we used the natural
identifications $H^0(\O(n)) \cong \V_n$ and $H^0(\O(n) \oplus \O(n))
\cong \V_{1,n}$ for any integer $n \geq 0$.

Let us define the vector subspaces $\V'$ and $\V''$ of $\V_{1,k}$ by \[
\V' := \{ x \otimes u_x + y \otimes u_y\ |\ u \in \V_{k+1} \}, \] where
$u_x$ and $u_y$ denote partial derivatives, and \[ \V'' :=
\und{\imath}(\V_{k-1}) = \{ x \otimes (y \cdot u) - y \otimes (x \cdot
u)\ |\ u \in \V_{k-1} \}.\] Then $\V_{1,k} = \V' \oplus \V''$ is easily
verified.

\proclaim{Proposition 3.4} Let $C \in \Y \subseteq \Z$ be a rational
contact curve and consider bundle isomorphisms $\phi_C, \phi'_C$ and
$\phi''_C$ which induce the commutative diagrams \tg{3-3} and \tg{3-4}.
Then we have \[ \matrix \phi_C^*(\V') = T_C\Y & \text{and} &
\phi_C^*(\V'') = \imath (H^0 (L^{\perp}_{|C} /\tau)). \endmatrix
\tag3-5$$
\endproclaim

\demo{Proof} The map $\eta: \V_{k+1} \to \V_{1,k}$ given by $\eta(u) :=
\frac 1 {k+1} (x \otimes u_x + y \otimes u_y)$ splits the top exact
sequence of \tg{3-4}. Thus, Corollary 3.3. implies that there is a map
$\delta: \V_{k+1} \to \V_{k-1}$ such that $T_C\Y = \phi_C^* \left( \{
(\eta + \und{\imath}\circ \delta)(u)\ |\ u \in \V_{k+1} \} \right)$.

Since $\V_{k+1}^* \otimes \V_{k-1} \cong \V_{2k} \oplus \V_{2k-2}
\oplus \ldots \oplus \V_2$ it follows that there are polynomials $v_i
\in \V_i$, $i = 2, 4, \ldots 2k$ such that \[ \delta(u) = \pair u
{v_{2k}} {k+1} + \pair u {v_{2k-2}} k + \ldots + \pair u {v_2} 2.\]

{}From Lemma 3.2. and an easy calculation we conclude that $\delta$ must
satisfy the following condition: \[ \text{if $r^2 | u$ for some $r \in
\V_1$, $u \in \V_{k+1}$ then $r | \delta(u)$}. \tag3-6$$

Using the $\Sl$\=equivariance of $\pair \_ \_ {}$ we compute that for
any $r \in \V_1$, $r | \delta(r^{k+1})$ if and only if $r | v_{2k}$.
Thus, \tg{3-6} implies that every $r \in \V_1$ divides $v_{2k}$, hence
$v_{2k} = 0$.

Next, a similar calculation shows that $r | \delta(r^k s)$ for all $s
\in \V_1$ if and only if $r|v_{2k-2}$. Again, this together with
\tg{3-6} implies $v_{2k-2} = 0$.

Continuing with similar arguments, we see successively that $v_{2k} =
v_{2k-2} = \ldots = v_2 = 0$, thus $\delta = 0$, and this shows the
first equation in \tg{3-5}. The second equation is immediate from the
commutativity of \tg{3-4}. \qed
\enddemo

\proclaim{Proposition 3.5} Let $F_{\Y}:= \pi^{-1} (\Y)$ with the
principal $G$-bundle $\pi: F \to \Z$ from \tg{2-1} and let $\pi_{\Y}:
\F_{\Y} \to \Y$ denote the total $\V_{k+1}$-coframe bundle of $\Y$. The
set \[ \hat F:= \left\{ \phi_C: \O(k) \oplus \O(k) \to
N_C\ \left|\ \matrix C \in \Y,\ \phi_C \text{\ a bundle
isomorphism}\\ \text{which satisfies \tg{3-5}}. \endmatrix \right.
\right\} \subseteq F_{\Y} \] is a reduction of $F_{\Y}$ with structure
group \[ G^{\Delta}:= \{(c \cdot A, A)\ |\ c \in \Bbb C^*, A \in \Sl \}
\subseteq G. \] Moreover, the map \[ \matrix \format \c\ & \c\ &
\l\\ \zeta: & \hat F & \to \F_{\Y}\\ & \phi_C & \mapsto (\und{pr} \circ
(\phi_C^*)^{-1}: T_C\Y \to \V_{k+1}) \endmatrix \] is an imbedding, and
the image $\zeta(\hat F) \subseteq \F_{\Y}$ is a $G_{k+1}$-structure on
$\Y$.
\endproclaim

\demo{Proof} The proof is straightforward: first of all, by our
previous discussion we know that $\pi^{-1}(C) \cap \hat{F} \neq
\emptyset$ for all $C \in \Y$. Moreover, if $\phi_C^1, \phi_C^2 \in
\hat F_C$, then $\psi:= (\phi_C^1)^{-1} \circ \phi_C^2 \in Aut(\O(k)
\oplus \O(k))$ must satisfy $\psi^*(\V') = \V'$ and $\psi^*(\V'') =
\V''$. This is the case precisely if $\psi \in G^{\Delta}$. Thus $\hat
F$ is a $G^{\Delta}$-reduction of $F_{\Y}$.

The verification of the stated properties of $\zeta$ is left to the
reader. \qed
\enddemo

By abuse of notation, we shall identify $\hat F$ with $\zeta(\hat F)$
and thus regard $\hat F$ as a $G_{k+1}$-structure on $\Y$. The
tautological 1-form of $\pi: \hat F \to \Y $ is then given by \[ \hat
\th = \und{pr} \circ (\th|_{\hat F}). \]

The decomposition $\V_{1,k} = \V' \oplus \V''$ induces the
decomposition \[ \V_{1,k}^* \otimes \V_{1,k} = ({\V'}^* \otimes \V')
\oplus ({\V'}^* \otimes \V'') \oplus ({\V''}^* \otimes \V') \oplus
({\V''}^* \otimes \V''). \]

Projection onto the first direct summand composed with conjugation by
$\und{pr}|_{\V'}$ yields a homomorphism \[ p: \gl{\V_{1,k}} \to
\gl{\V_{k+1}}. \]

It is not hard to verify that the 1-form \[ \hat{\om} := (p \circ
\om)|_{\hat F} \] yields a connection on $\pi: \hat F \to \Y$.

\definition{Definition 3.6}
\roster
\item A $G_{1,k}$-connection is a triple $(\pi: F \to \Z, \th, \om)$ of
a $G_{1,k}$-reduction $F$ of $\Z$ and the tautological and connection
1-forms $\th$ and $\om$.

\item Likewise, a $G_{k+1}$-connection is a triple $(\pi: \hat F \to
\Y, \hat \th, \hat \om)$ of a $G_{k+1}$-reduction $\hat F$ of $\Y$ and
the tautological and connection 1-forms $\hat \th$ and $\hat \om$.

\item Suppose $(\pi: F \to \Z, \th, \om)$ and $(\pi: \hat F \to \Y,
\hat \th, \hat \om)$ are a $G_{1,k}$-connection and a
$G_{k+1}$-connection respectively and suppose there is an injective
bundle map \[ \CD \hat F @>\jmath>> F\\ @V{\pi}VV @V{\pi}VV\\ \Y @>\und
\jmath>> \Z \endCD \] such that $\hat \th = \und{pr} \circ \th|_{\hat
F}$ and $\hat \om = p \circ \om|_{\hat F}$.

Then $(\pi: \hat F \to \Y, \hat \th, \hat \om)$ is called a {\it
restriction} of $(\pi: F \to \Z, \th, \om)$, whereas $(\pi: F \to \Z,
\th, \om)$ is called an {\it extension} of $(\pi: \hat F \to \Y, \hat
\th, \hat \om)$.
\endroster
\enddefinition

The following Proposition is straightforward and the proof is omitted.

\proclaim{Proposition 3.7} If the $G_{1,k}$-connection $(\pi: F \to \Z,
\th, \om)$ is an extension of the $G_{k+1}$-connection $(\hat \pi: \hat
F \to \Y, \hat \th, \hat \om)$ and if $\Theta$ and $\hat{\Theta}$
denote the torsion of $\om$ and $\hat{\om}$ respectively then \[
\hat{\Theta} = \und{pr} \circ \Theta|_{\hat F}. \text{\qed} \]
\endproclaim

\definition{Definition 3.8} A connection on $\pi_{\Y}: \hat F \to \Y$
is called {\it special} if it is the restriction of a special
connection on $\pi: F \to \Z$.
\enddefinition

Of course, from our discussion preceding Definition 3.6. we know that
if $\Z$ is the moduli space of rational curves in $\W$ whose normal
bundle is equivalent to $\O(k) \oplus \O(k)$ and if $\Y \subseteq \Z$
is the subset of contact curves then every connection on $\pi: F \to
\Z$ has a restriction to the $G_{k+1}$-structure $\pi: \hat F \to \Y$.

Let us now investigate the {\it intrinsic torsion} of both the
$G_{1,k}$-structure $\pi: F \to \Z$ and the $G_{k+1}$-structure
$\pi_{\Y}: \hat F \to \Y$. The Spencer sequence \tg{1-1} reads \[
\matrix & 0 \longto \g{1,k}^{(1)} \longto \V_{1,k}^* \otimes \g{1,k}
\ltop{Sp} \Lambda^2\V_{1,k}^* \otimes \V_{1,k} \longto H^{0,2}(\g{1,k})
\longto 0\\ \text{and}\\ & 0 \longto \g {k+1}^{(1)} \longto \V_{k+1}^*
\otimes \g {k+1} \ltop{Sp} \Lambda^2\V_{k+1}^* \otimes \V_{k+1} \longto
H^{0,2}(\g {k+1}) \longto 0. \endmatrix \tag3-7$$

\proclaim{Lemma 3.9} If $k \geq 2$ then $\g{1,k}^{(1)} = 0$ and $\g
{k+1}^{(1)} = 0$.
\endproclaim

\demo{Proof} Let $\varphi \in \g{1,k}^{(1)}$. We regard $\varphi$ as a
linear map $\varphi: \V_{1,k} \to \g{1,k}$. Pick two arbitrary bases
$(r_1, r_2)$ and $(s_1, s_2)$ of $\V_1$. Then the set $\{ r_i \otimes
s_1^{k-j} s_2^j\ |\ i = 1,2, j = 0, \ldots k \}$ forms a basis of
$\V_{1,k}$. We have \[ \left(\varphi(r_1 \otimes s_1^k) \right) (r_2
\otimes s_2^k) - \left(\varphi(r_2 \otimes s_2^k) \right) (r_1 \otimes
s_1^k) = 0. \]

If we let $\varphi(r_i \otimes s_i^k) := (A_i, B_i)$ for $i=1,2$ be the
decomposition in $\g{1,k} \cong \glC2 \oplus \slC$, then this equation
reads \[ ((A_1 \cdot r_2) \otimes s_2^k) + (r_2 \otimes (B_1 \cdot
s_2^k)) - ((A_2 \cdot r_1) \otimes s_1^k) - (r_1 \otimes (B_2 \cdot
s_1^k)) = 0. \tag3-8$$

Note that $B_1 \cdot s_2^k \in span\{ s_1 s_2^{k-1}, s_2^k\}$. Taking
the $(r_2 \otimes s_1 s_2^{k-1})$-component of \tg{3-8} w.r.t. the
above basis we conlude that $s_2^k$ is an eigenvector of $B_1$. Since
this is true for {\it any} $s_2$ which is linearly independent of
$s_1$, it follows that $B_1$ is a multiple of the identity. On the
other hand, $trace(B_1) = 0$. Thus, we have $B_1 = 0$.

Likewise, $A_1 = 0$, hence $\varphi(r_1 \otimes s_1^k) = 0$ for
arbitrary $r_1, s_1 \in \V_1$. Since elements of this form span all of
$\V_{1,k}$, $\varphi = 0$ follows.

The proof of the second statement is of similar nature but simpler. We
omit the details. \qed
\enddemo

To calculate the irreducible components of \tg{3-7}, note that as a
$G$-module, $\g{1,k} \cong \V_{0,0} \oplus \V_{2,0} \oplus \V_{0,2}$.
In fact, the equivalence is determined by the equation \[ (p_{0,0} +
p_{2,0} + p_{0,2}) \cdot q_{1,k} := \pair {p_{0,0}} {q_{1,k}} {0,0} +
\pair {p_{2,0}} {q_{1,k}} {1,0} + \pair {p_{0,2}} {q_{1,k}} {0,1}
\tag3-9$$ for all $p_{i,j} \in \V_{i,j}$ and $q_{1,k} \in \V_{1,k}$.

For the rest of this section we shall assume that $k=2$. In this case,
a calculation shows that the decomposition of the Spencer sequence
\tg{3-7} into irreducible submodules is \[ 0 \longto \matrix \V_{1,0}
\oplus 3 \V_{1,2} \oplus \V_{1,4}\\ \oplus \V_{3,2} \endmatrix
\ltop{Sp} \matrix \V_{1,0} \oplus 3 \V_{1,2} \oplus 2 \V_{1,4} \oplus
\V_{1,6}\\ \oplus \V_{3,0} \oplus \V_{3,2} \oplus \V_{3,4} \endmatrix
\longto \matrix \V_{1,4} \oplus \V_{1,6}\\ \oplus \V_{3,0} \oplus
\V_{3,4} \endmatrix \longto 0. \]

More explicitly, if $\varphi \in \V_{1,2}^* \otimes \g{1,2}$, then
there are elements $r_{i,j}, r'_{i,j}, r''_{i,j} \in \V_{i,j}$ such
that for $p_{1,2} \in \V_{1,2}$,

\[ \aligned \varphi(p_{1,2}) = & \left(\pair {r_{1,2}} {p_{1,2}}
{1,2}\right) + \left(\pair {r_{3,2}} {p_{1,2}} {1,2} + \pair {r'_{1,2}}
{p_{1,2}} {0,2} \right)\\ & \ + \left(\pair {r_{1,4}} {p_{1,2}} {1,2} +
\pair {r''_{1,2}} {p_{1,2}} {1,1} + \pair {r_{1,0}} {p_{1,2}}
{1,0}\right) \in \V_{0,0} \oplus \V_{2,0} \oplus \V_{0,2}. \endaligned
\tag3-10$$

Likewise, for any $T \in \Lambda^2\V_{1,2}^* \otimes \V_{1,2}$, there
are elements $s_{i,j}, s'_{i,j}, s''_{i,j} \in \V_{i,j}$ such that for
all $p,q \in \V_{1,2}$,

\[ \aligned T(p, q) = & \pair {s_{1,2}} {\pair p q {1,0}} {0,2} + \pair
{s_{1,4}} {\pair p q {1,0}} {0,3} + \pair {s_{1,6}} {\pair p q {1,0}}
{0,4}\\ & \qquad + \pair {s_{1,0}} {\pair p q {0,1}} {1,0} + \pair
{s'_{1,2}} {\pair p q {0,1}} {1,1} + \pair {s'_{1,4}} {\pair p q {0,1}}
{1,2}\\ & \qquad + \pair {s_{3,0}} {\pair p q {0,1}} {2,0} + \pair
{s_{3,2}} {\pair p q {0,1}} {2,1} + \pair {s_{3,4}} {\pair p q {0,1}}
{2,2}\\ & \qquad + \pair {s''_{1,2}} {\pair p q {1,2}} {0,0}.
\endaligned \tag3-11$$

Using the tuple $(s_{1,2}, s_{1,4}, s_{1,6}, s_{1,0}, s'_{1,2},
s'_{1,4}, s_{3,0}, s_{3,2}, s_{3,4}, s''_{1,2})$ as coordinates for
$\Lambda^2\V_{1,2}^* \otimes \V_{1,2}$, another calculation shows that
\[ \aligned Sp(\varphi) = & \Big( -\frac16 (r_{1,2} - 3 r'_{1,2} + 4
r''_{1,2}),\ - \frac12 r_{1,4},\ 0,\ r_{1,0},\ -\frac18 (r_{1,2} +
r'_{1,2} - 4 r''_{1,2}),\\ & \qquad - \frac12 r_{1,4},\ 0,\ - \frac14
r_{3,2},\ 0,\ - \frac13 (r_{1,2} - 3 r'_{1,2} -8 r''_{1,2}) \Big),
\endaligned \tag3-12$$ where the $r_{i,j}$'s are determined by
$\varphi$ as in \tg{3-10}.

\proclaim{Lemma 3.10} Let $\Z$ be as before, and suppose that $\om_1$
is a special connection on $\Z$. Then another connection $\om_2$ on
$\Z$ is special if and only if there are functions $r_{i,j}, r_{i,j}'$
on $F$ with values in $\V_{i,j}$ such that \[ \om_2 = \om_1 +
\left(\pair {r_{1,2}} {\theta} {1,2}\right) + \left(\pair {r_{3,2}}
{\theta} {1,2} + \pair {r'_{1,2}} {\theta} {0,2} \right) + \left(\pair
{r_{1,0}} {\theta} {1,0}\right). \tag3-13$$ Here, we use the
identification \tg{3-9} to regard the $\om_i$'s as $\V_{0,0} \oplus
\V_{2,0} \oplus \V_{0,2}$-valued 1-forms on $F$.
\endproclaim

\demo{Proof} First of all, we define the vector bundles $\V_{i,j}^F :=
F \times_G \V_{i,j}$ over $\Z$. Recall from \S1 that the difference
between two connections on $F \to \Z$ is determined by a section of the
bundle $T^*\Z \otimes \g F$. Since by \tg{3-9} we have $\g F \cong
\V_{0,0}^F \oplus \V_{2,0}^F \oplus \V_{0,2}^F$, we can decompose \[
\matrix T^*\Z \otimes \g F = B_1^F \oplus B_2^F & \text{with} & \left\{
\matrix B_1^F := T^*\Z \otimes \left(\V_{0,0}^F \oplus
\V_{2,0}^F\right), & \text{and}\\ B_2^F := T^*\Z \otimes \V_{0,2}^F.
\endmatrix \right. \endmatrix \]

Note that $B_2^F \cong \V_{1,2}^F \otimes \V_{0,2}^F \cong \V_{1,0}^F
\oplus \V_{1,2}^F \oplus \V_{1,4}^F$.

Consider the vector bundle \[ \Delta := \bigcup_{t \in \Z}
H^0(Hom(N_{C_t}, \tau_{C_t})) \to \Z. \]

For a fixed $t \in \Z$, $Hom(N_{C_t}, \tau_{C_t}) \cong Hom(\O(2)
\oplus \O(2), \O(2)) \cong \O(0) \oplus \O(0)$, hence $H^0(Hom(N_{C_t},
\tau_{C_t})) \cong \V_{1,0}$ as a $G$-module, and thus $\Delta \cong
\V_{1,0}^F$.

Let us fix a special connection $\om_0 = \tilde{\om} + \hat{\sigma}_0$
with some section $\sigma_0$ of the split bundle $S \to \Z$. Given a
local section $\delta$ of $\Delta$, we let $\om := \tilde{\om} +
\hat{\sigma}$ with $\sigma := \sigma_0 + \delta$, and define
$\psi(\delta) := \om - \om_0$. From Lemma 2.4. it is easy to verify
that $\psi(\delta)$ is a local section of $B_2^F \subseteq T^*\Z
\otimes \g F$, that the correspondence $\delta \mapsto \psi(\delta)$
determines a bundle map $\psi: \Delta \to B_2^F$, and that $\psi$ is
independent of the choice of $\om_0$. Also, it is obvious that $\psi$
is non-vanishing, hence by Schur's Lemma $\psi(\Delta) = \V_{1,0}^F
\subseteq B_2^F$.

Let $\om_1 = \tilde{\om}_1 + \hat{\sigma}_1$ be the decomposition of
the special connection $\om_1$ where $\sigma_1$ is a section of the
split bundle $S$. Then $\om_2$ is {\it special} if and only if $\om_2 =
\tilde{\om}_2 + \hat{\sigma}_2$ for some section $\sigma_2$ of $S$, if
and only if $\om_2 - \om_1 = \left(\tilde{\om}_2 - \tilde{\om}_1\right)
+ \psi(\delta)$ where $\delta := \sigma_2 - \sigma_1$ is a section of
$\Delta$, if and only if $\om_2 - \om_1$ is a section of $B_1^F \oplus
\psi(\Delta) = B_1^F \oplus \V_{1,0}^F \subseteq T^*\Z \otimes \g F$.

Comparing \tg{3-13} with \tg{3-10}, we see that this is satisfied if
and only if \tg{3-13} holds. \qed
\enddemo

\proclaim{Theorem 3.11} Let $\Z$ be the moduli space of rational curves
in a 3-fold $\W$ whose normal bundle is equivalent to $\O(2) \oplus
\O(2)$, and let $\pi: F \to \Z$ be the associated $G_{1,2}$-structure.
Then there is a unique connection $\om$ on $F$ and a function $s_{3,0}:
F \to \V_{3,0}$ such that the torsion of $\om$ is given by \[ \Theta =
\pair {s_{3,0}(u)} {\pair {\theta} {\theta} {0,1}} {2,0}. \]

Moreover, $\om$ is a special connection.
\endproclaim

This simple form of the torsion is quite remarkable; indeed, Lemma 3.9.
implies that for $k=2$, $rank(H^{0,2}_F) = 48$. Thus, Theorem 3.11.
says that {\it most} of the intrinsic torsion of $F$ vanishes.

\demo{Proof} First of all, note that it suffices to prove the Theorem
locally, i.e. we need to show that $F$ can be covered by open sets
$\U_i$ on which a connection $\om_i$ with the stated properties exists.
Then, by uniqueness, $\om_i$ and $\om_j$ must coincide on $\U_i \cap
\U_j$, thus $\om_i$ is the restriction to $\U_i$ of a connection $\om$
defined on {\it all} of $F$. In the proof, we will replace $\U_i$ by
$\Z$ and thus we may assume that all local properties of $\Z$ hold
globally.

By Proposition 2.6. we can find a {\it special} connection $\om_0$ on
$F$. Then there are functions $(s_{1,2}, s_{1,4}, \ldots, s''_{1,2})$
on $F$ with values in $\V_{1,2}, \V_{1,4}, \ldots, \V_{1,2}$
respectively such that the torsion $\Theta_0$ of $\om_0$ is given by
$\Theta_0 (\xi, \xi') = T_0(p,q)$ with $p = \theta(\xi), q =
\theta(\xi')$, and where $T_0(p,q)$ is determined by the $s_{i,j}$'s as
in \tg{3-11}.

We call $T_0: F \to \Lambda^2\V_{1,2}^* \otimes \V_{1,2}$ the {\it
torsion map of $\om_0$}. By abuse of notation, we write $T_0 = s_{1,2}
+ s_{1,4} + \ldots + s''_{1,2}$.

We shall say that $r \in \V_1$ {\it divides} $p \in \V_{1,k}$ and write
$r|p$ if $p = x \otimes p_1 + y \otimes p_2$ and $r$ divides both $p_1$
and $p_2$. From Corollary 2.8. we obtain the following criterion: \[
\text{If for $r \in \V_1$ and $p,q \in \V_{1,2}$ we have $r|p$ and
$r|q$ then also $r|T_0(u)(p,q)$.} \tag3-14$$

{}From \tg{3-14}, we can conclude that some of the $s_{i,j}$'s must
vanish identically. For example, let $p := x \otimes r^2$ and $q := y
\otimes r^2$ for $r \in \V_1$. Then $\pair p q {0,1} = \pair p q {1,2}
= 0$ and $\pair p q {1,0} = 1 \otimes r^4 \in \V_{0,4}$. Since $r|p$
and $r|q$, \tg{3-14} implies that $r|T_0(u)(p,q)$. An easy computation
shows that this is the case if and only if $r|s_{1,6}(u)$. But this
must be true for {\it all} $r \in \V_1$ and $u \in F$, thus we conclude
$s_{1,6} = 0$.

By similar calculations we see that \tg{3-14} is satisfied if and only
if \[ s_{1,4} = s'_{1,4} = s_{1,6} = s_{3,4} = s''_{1,2} - 2 s_{1,2} =
0. \tag3-15$$

{}From here it follows that the intrinsic torsion of the
$G_{1,2}$-structure is represented by $s_{3,0}$. Also, from \tg{3-12}
and \tg{3-15} it follows that $T_0 - s_{3,0} = Sp(\varphi)$ for some
function $\varphi: F \to \V_{1,2}^* \otimes \g{1,2}$ whose $r_{1,4}$
and $r''_{1,2}$ component vanish. Then Lemma 3.10. implies that the
connection $\om := \om_0 - \varphi$ is still {\it special}, and if we
let $T$ be the torsion map of $\om$ then by \tg{1-7} we have $T = T_0 -
Sp(\varphi) = s_{0,3}$. Thus, the torsion of the special connection
$\om$ is of the desired form.

The uniqueness follows from Proposition 1.2. together with Lemma 3.9.
\qed
\enddemo

\definition{Definition 3.12} The unique special connection from Theorem
3.11. is called the {\it intrinsic connection of $\Z$}.
\enddefinition

\proclaim{Theorem 3.13} Let $\W$ be a contact 3-fold with contact line
bundle $L \to \W$, and let $\Y$ be the moduli space of rational contact
curves $C$ such that $L|_C \cong \O(-3)$. Then there is a unique
torsion free connection $\hat{\om}$ on the $G_3$-structure $\pi_{\Y}:
\hat F \to \Y$. Moreover, $\hat{\om}$ is special.
\endproclaim

\demo{Proof} By Proposition 3.1. $\Y \subseteq \Z$ is a submanifold
with $\Z$ as in Theorem 3.11. Let $\om_0$ be the intrinsic connection
on $\Z$, and let $s_3: F \to \V_3$ be such that $s_{3} \otimes 1: F \to
\V_{3,0}$ is the torsion function of $\om_0$ from Theorem 3.11.

Let $s_{1,2}:= x \otimes (s_3)_x + y \otimes (s_3)_y \in \V_{1,2}$, and
let $\varphi: F \to \V_{1,2}^* \otimes \g{1,2}$ be determined by
\tg{3-10} with $r_{1,2} = r'_{1,2} := 2 s_{1,2}$, all other $r_{i,j}$'s
$=0$. As before, $\varphi$ can be regarded as a section of $T^*\Z
\otimes \g F$, and by Lemma 3.10, the connection $\om := \om_0 +
\varphi$ is again a special connection on $\Z$. We denote by
$\hat{\om}$ the restriction of $\om$ to $\Y$, and let $\Theta,
\Theta_0$ and $\hat{\Theta}$ denote the torsion forms of $\om, \om_0$
and $\hat{\om}$ respectively. By \tg{1-7}, $\Theta = \Theta_0 +
Sp(\varphi)$.

Then by Proposition 3.7. and some calculation we have

\[ \aligned \hat{\Theta} & = \und{pr} \circ \Theta\\ & = \und{pr}
\pmatrix \pair {s_{3,0}} {\pair{\theta} {\theta} {0,1}} {2,0} - \frac12
\pair {s_{1,2}} {\pair{\theta} {\theta} {0,1}} {1,1}\\ \qquad + \frac23
\pair {s_{1,2}} {\pair{\theta} {\theta} {1,0}} {0,2} + \frac43 \pair
{s_{1,2}} {\pair{\theta} {\theta} {1,2}} {0,0}\endpmatrix\\ & = 0.
\endaligned \]

Thus, $\hat{\om}$ is the desired torsion free special connection. The
uniqueness follows from Proposition 1.2. together with Lemma 3.9. \qed
\enddemo

\subheading{\S4 Torsion free $G_{1,2}$-structures}

In this entire section, we shall consider complex sixfolds $\Z$ which
carry a {\it torsion free} $G_{1,2}$-structure $\pi: F \to \Z$. In this
case, there is a unique torsion free connection $\om = \om_{0,0} +
\om_{2,0} + \om_{0,2}$ on $F$ where $\om_{i,j}$ takes values in
$\V_{i,j}$. Here, we used the identification $\g{1,2} \cong \V_{0,0}
\oplus \V_{2,0} \oplus \V_{0,2}$ from \tg{3-9}.

For convenience, we shall define the pairings \[ \matrix \ppair {\_}
{\_} k: & (\V_{0,0} \oplus \V_{2,0} \oplus \V_{0,2}) \otimes \V_{i,j} &
\to & \V_{i,j}\\ & \ppair {p_{0,0} + p_{2,0} + p_{0,2}} {\ q} k & := &
k \pair {p_{0,0}} q {0,0} + \pair {p_{2,0}} q {1,0} + \pair {p_{0,2}} q
{0,1} \endmatrix \]

Then the {\it first structure equation} of $\om$ reads \[ d\theta +
\ppair {\om} {\theta} 1 = 0 \tag4-1$$ with the $\V_{1,2}$-valued
tautological 1-form $\theta$.

Moreover, the {\it curvature 2-form $\Omega$} takes values in $\g{1,2}
\cong \V_{0,0} \oplus \V_{2,0} \oplus \V_{0,2}$, and is defined as \[
\aligned \Omega & = d\om + \om \w \om\\ & = d\om - \frac12 \left( \pair
{\om_{2,0}} {\om_{2,0}} {1,0} + \pair {\om_{0,2}} {\om_{0,2}} {0,1}
\right). \endaligned \tag4-2$$

Differentiating \tg{4-1}, we obtain the {\it first Bianchi identity} \[
\ppair {\Omega} {\theta} 1 = 0. \tag4-3$$

Let $\K (\g{1,2})$ be given by the exact sequence \[ 0 \longto \K
(\g{1,2}) \longto \Lambda^2 \V_{1,2}^* \otimes \g{1,2} \ltop{Sp_2}
\Lambda^3 \V_{1,2}^* \otimes \V_{1,2}, \] where $Sp_2$ is given by
skew-symmetrization of $\Lambda^2 \V_{1,2} \otimes \g{1,2} \subseteq
\Lambda^2 \V_{1,2} \otimes (\V_{1,2}^* \otimes \V_{1,2})$. The first
Bianchi identity \tg{4-3} can be interpreted as stating that $\Omega$
is a section of $F \times_G \K (\g{1,2})$.

A calculation shows that, as a $G$-module, $\K (\g{1,2}) \cong \V_{2,0}
\oplus \V_{0,2}$. More explicitly, there is a function $\a = a_{2,0} +
a_{0,2}: F \to \V_{2,0} \oplus \V_{0,2}$ such that \[ \aligned \Omega &
= \left( -4 \pair {a_{2,0}} {\pair {\theta} {\theta} {1,2}} {0,0} + 3
\pair {a_{0,2}} {\pair {\theta} {\theta} {0,1}} {0,2} \right) \\ &
\qquad + \left( \pair {a_{2,0}} {\pair {\theta} {\theta} {0,1}} {2,0} +
\pair {a_{0,2}} {\pair {\theta} {\theta} {1,0}} {0,2} -7 \pair
{a_{0,2}} {\pair {\theta} {\theta} {1,2}} {0,0} \right). \endaligned
\tag4-4$$

This implies, in particular, that $d\om_{0,0} = 0$. Therefore, by the
{\it Ambrose-Singer-Holonomy Theorem} \cite{KN}, the holonomy of $\om$
is contained in the subgroup \[ H_{1,2} := G_{1,2} \cap {\text{\it
Sl}(V_{1,2})}. \]

Taking the derivative of \tg{4-4} and solving for $d\a$, we see that
there is a function $\bb: F \to \V_{1,2}$ such that \[ d\a = \ppair
{\om} \a {-2} + 3 \pair {\bb} {\theta} {0,2} + \pair {\bb} {\theta}
{1,1}. \tag4-5$$

Once again, we take the derivative of \tg{4-5} and solve for $d\bb$. We
see that there is a function $c: F \to V_{0,0}$ such that \[ \aligned
d\bb & = \ppair {\bb} {\om} {-3} + 2 \pair {\pair {a_{2,0}} {a_{0,2}}
{0,0}} {\th} {1,1} + \pair {\pair {a_{0,2}} {a_{0,2}} {0,0}} {\th}
{0,2}\\ & \qquad + \pair {-\frac43 \pair {a_{2,0}} {a_{2,0}} {2,0} - 7
\pair {a_{0,2}} {a_{0,2}} {0,2} + c} {\th} {0,0} \endaligned \tag4-6$$

Taking exterior derivatives one more time and solving for $dc$ we
calculate that \[ dc = -4 c \om_{0,0}. \tag4-7$$

The reader who is familiar with \cite{Br} will note the similarity of
the structure equations \tg{4-1} - \tg{4-7} with the structure
equations for $H_3$-connections where $H_3 = G_3 \cap {\text{\it
Sl}(\V_3)}$. This is by no means a coincidence. As we shall see in the
following section, there is a close relationship between
$H_{1,2}$-structures and $H_3$-structures.

Let $F_0 \subseteq F$ be an integral hypersurface of $\om_{0,0}$, i.e.
a hypersurface such that $\om_{0,0}|_{F_0} \equiv 0$. Then $F_0$ is a
torsion free $H_{1,2}$-reduction of $F$. We shall denote the
restrictions of $\theta$, $\om_{2,0}$, $\om_{0,2}$, $\a$ and $\bb$ to
$F_0$ by the same letters. Note that $\om_0 := \om|_{F_0} = \om_{2,0} +
\om_{0,2}$ is $\V_{2,0} \oplus \V_{0,2}$-valued. Also, by \tg{4-7}, $c$
is {\it constant} on $F_0$.

\definition{Definition 4.1} Let $\pi: F \to \Z$ be a torsion free
$G_{1,2}$-structure and let $F_0 \subseteq F$ be an integral
hypersurface of $\om_{0,0}$. Then $F_0$ is called an {\it associated
$H_{1,2}$-structure of $F$}.
\enddefinition

The choice of associated $H_{1,2}$-structures is, of course, not
unique. However, given two such structures $F_0$ and $F_0'$ then $F_0'
= R_{t I} \cdot F_0$ for some $t \in \Bbb C^*$. {\it Hence, to each
torsion free $G_{1,2}$-structure there is a one parameter family of
associated $H_{1,2}$-structures}.

Our approach to solve the structure equations \tg{4-1} - \tg{4-7} will
be motivated by the steps pursued in \cite{Br} to solve the structure
equations of an $H_3$-connection.

Let \[ \matrix K := \a + \bb: F_0 \to \V, & \text{where} & \V =
\V_{2,0} \oplus \V_{0,2} \oplus \V_{1,2}. \endmatrix \]

Equations \tg{4-5} - \tg{4-6} can be summarized as \[ dK = J (\theta +
\om_0) \] where $J$ is a function on $F$ with values in $Hom(\V,\V)$.
Now $J = K^*(\J_c)$ where $\J_c: \V \to Hom(\V,\V)$ is a polynomial
mapping which depends upon a parameter $c$. If we write $\J_c$ relative
to the standard basis of $\V$ then it has a $12 \times 12$-matrix
representation whose entries are polynomials in the components of $\a$
and $\bb$.

It turns out that this matrix $\J_c$ is not invertible. In fact,
generically the rank of $\J_c$ is calculated to be 10. This implies
that the image of $K$ is contained in some 10-dimensional subvariety of
$\V$.

Let us this once comment on the mechanical calculations which are
performed to arrive at this conclusion. The attempt of simply taking
the determinant of $\J_c$ on \math failed miserably at first: after
more than 10 minutes of calculation, memory overflows occured.

The next approach was to use the $H_{1,2}$-equivariance of $\J_c$.
Under the generic assumption that both $a_{2,0}$ and $a_{0,2}$ are not
squares of a linear polynomial, we may assume that $a_{2,0} = t x y
\otimes 1$ and $a_{0,2} = t' x y \otimes 1$ for some $t, t' \in \Bbb
C$. Making this replacement simplifies $\J_c$ to a matrix $\J_c'$ of
equal rank which is drastically simpler, and calculating that
$det(\J_c') = 0$ on \math is a matter of less than a minute.

Moreover, we can explicitly compute the kernel of $\J_c'$, and thus by
equivariance the kernel of $\J_c$. The result can be described as
follows. Let
\[ \matrix \matrix \matrix \format \r\ & \c\ & \l\\
d_1 \otimes 1 & := & \pair {a_{2,0}} {a_{2,0}} {2,0}\\
d_2 \otimes 1 & := & \pair {a_{0,2}} {a_{0,2}} {0,2} \endmatrix &
\matrix \format \r\ & \c\ & \l\\
e_1 \otimes 1 & := & \pair {\pair {a_{2,0}} {\bb} {1,0}} {\bb} {1,2}\\
e_2 \otimes 1 & := & \pair {\pair {a_{0,2}} {\bb} {0,1}} {\bb} {1,2}
\endmatrix &
\matrix \format \r\ & \c\ & \l\\
b_{0,2} & := & \pair {\bb} {\bb} {1,1}\\
b_{2,0} & := & \pair {\bb} {\bb} {0,2}\\
b_{2,4} & := & \pair {\bb} {\bb} {0,0}
\endmatrix \endmatrix \\
\matrix \format \r\ & \c\ & \l\\
p_{2,0} & := & \pair {\pair {a_{2,0}} {a_{0,2}} {0,0}} {a_{0,2}}
{0,2}\\
p_{2,4} & := & \pair {\pair {a_{2,0}} {a_{2,0}} {0,0}} {a_{2,0}}
{2,0}\\
p_{0,2} & := & 16 \pair {\pair {a_{2,0}} {a_{0,2}} {0,0}} {a_{2,0}}
{2,0}
+ 9 \pair {\pair {a_{0,2}} {a_{0,2}} {0,0}} {a_{0,2}} {0,2}\\
& & \qquad - 12\ c\ a_{0,2} + 3 b_{0,2}
\endmatrix
\endmatrix \]

Here, $d_i$ and $e_i$ are $\Bbb C$-valued functions on $\V$, while
$b_{i,j}$ and $p_{i,j}$ are functions on $\V$ with values in
$\V_{i,j}$. Then we have the

\proclaim{Proposition 4.2} Let \[ f_1^c := (4 d_1 - 9 d_2) (4 d_1 + 27
d_2 - 6 c) + 72 e_1 - 54 e_2\] and \[ f_2^c := 4 d_2 (4 d_1 + 9 d_2 - 3
c)^2 + 96 \pair {p_{2,0}} {b_{2,0}} {2,0} + 3 \pair {p_{0,2}} {b_{0,2}}
{0,2} + 48 \pair {p_{2,4}} {b_{2,4}} {2,4}.\] Then $d(K \circ f_i^c) =
0$ for $i = 1,2$, and hence $K$ maps $F_0$ into a level set of
$(f_1^c,f_2^c)$. Moreover, $rank(\J_c) = 10$ at $x \in \V$ if and only
if $df_1^c \w df_2^c|_x \neq 0$.
\endproclaim

\demo{Proof} The calculations involved to verify this Proposition were
all performed on \math and will not be presented here in further
detail. \qed
\enddemo

We let \[ \Sigma_c := \{ x \in \V\ |\ df_1^c \w df_2^c|_x = 0 \}. \]
Then by Proposition 4.2. we know that $rank(K)_u = 10$ at $u \in F_0$
if and only if $K(u) \notin \Sigma_c$.

Let us define the functions $r_{i,j}^k: \V \to \V_{i,j}$ for $k = 1,2$
by the equation \[ \matrix df_k^c = \frac16 \pair {r_{2,0}^k}
{da_{2,0}} {2,0} + \frac12 \pair {r_{0,2}^k} {da_{0,2}} {0,2} + \frac12
\pair {r_{1,2}^k} {d\bb} {1,2} & \text{for $k = 1,2$}, \endmatrix \]
and define the vector fields $Z_k$, $k = 1,2$, on $F_0$ by \[ \matrix
\om_0(Z_k) = r_{2,0}^k + r_{0,2}^k & \text{and} & \th(Z_k) = r_{1,2}^k.
\endmatrix \]

Then another \math calculation yields

\proclaim{Proposition 4.3} The vector fileds $Z_1$ and $Z_2$ on $F_0$
defined above are symmetries, i.e. their Lie derivatives satisfy \[
\matrix \frak L_{Z_k} (\om_0) = \om_0 & \text{and} & \frak L_{Z_k}
(\th) = \th & \text{for $k = 1,2$}. \endmatrix \tag4-8$$ Moreover,
$[Z_1, Z_2] = 0$. \qed
\endproclaim

\proclaim{Corollary 4.4} Either $rank(K) \equiv 10$ or $rank(K) < 10$
on all of $F_0$.
\endproclaim

\demo{Proof} From \tg{4-8}, standard arguments show that a symmetry
either vanishes {\it everywhere} or {\it nowhere} on $F_0$. Thus,
either $Z_1$ and $Z_2$ are pointwise linearly independent {\it
everywhere} or {\it nowhere} on $F_0$.

{}From the definitions of the $Z_k$'s it follows that $Z_1$ and $Z_2$ are
linearly independent if and only if $df_1^c$ and $df_2^c$ are linearly
independent. The claim follows then from Proposition 4.2. \qed
\enddemo

\definition{Definition 4.5} A torsion free $H_{1,2}$-connection is
called {\it regular} if $rank(K) \equiv 10$, with the map $K: F_0 \to
\V$ from above.

A torsion free $G_{1,2}$-connection is called {\it regular} if one and
hence all of its associated $H_{1,2}$-structures are regular.
\enddefinition

Thus, for a regular $H_{1,2}$-connection, the map $K$ is a submersion
onto an open subset of the regular part of a level set of $(f_1^c,
f_2^c)$ in $\V$.

\definition{Definition 4.6} Given constants $c, c_1, c_2 \in \Bbb C$
let \[ \C(c,c_1,c_2) := (f_1^c,f_2^c)^{-1} (c_1,c_2) \backslash
\Sigma_c \subseteq \V. \]

If for a regular torsion free $H_{1,2}$-connection on $\pi: F_0 \to \Z$
the image of $K: F_0 \to \V$ is contained in $\C_{c, c_1, c_2}$ then we
call the triple $(c, c_1, c_2)$ the {\it structure constants of the
connection}.
\enddefinition

Let us now consider the question of {\it existence} of torsion free
$H_{1,2}$-connections.

\proclaim{Theorem 4.7} Given constants $c, c_1, c_2 \in \Bbb C$ let $\C
:= \C(c,c_1,c_2) \subseteq \V$. Then $\C$ can be covered by open sets
$U$ which have the following property: there exists a holomorphic
principal $\Bbb C^2$-bundle $K: F_0 \to U$ over $U$ and holomorphic
1-forms $\th$ and $\om_0$ on $F_0$ with values in $\V_{1,2}$ and
$\V_{2,0} \oplus \V_{0,2}$ respectively satisfying

\roster
\item the $\V_{1,2} \oplus \V_{2,0} \oplus \V_{0,2}$-valued 1-form $\th
+ \om_0$ is a coframe on $F_0$,
\item equations \tg{4-1} - \tg{4-7} are satisfied if we set $\om_{0,0}
= 0$, and if $K = \a + \bb$ is the decomposition of $K$ into its
$\V_{2,0} \oplus \V_{0,2}$ and $\V_{1,2}$ components respectively.
\endroster

Moreover, the triple $(F_0, \th, \om_0)$ is unique in the sense that if
$(F_0', \th', \om_0')$ is another triple satisfying \tg1 and \tg2, then
there is a bundle isomorphism between $F_0$ and $F_0'$ which identifies
the coframings.
\endproclaim

\demo{Proof} Let $\ov a, \ov b, {\ov r}_{i,j}^k$ and $\ov J$ be the
restrictions of the functions $\a, \bb, r_{i,j}^k$ and $\J_c$
respectively to $\C$.

By definition of $\C$ we have $rank(\ov J) \equiv 10$. From here it
follows that there exist smooth 1-forms $\ov \th$ and ${\ov \om}_0$ on
$\C$ with values in $\V_{1,2}$ and $\V_{2,0} \oplus \V_{0,2}$
respectively such that \[ d\ov a + d\ov b = \ov J(\ov \th + {\ov
\om}_0). \tag4-9$$

Since $d\ov a, d\ov b$ and $\ov J$ are holomorphic, we may assume that
$\ov \th$ and $\ov \om_0$ are of type $(1,0)$.

The kernel of $\ov J$ is spanned at each point by the vectors ${\ov
r}_{2,0}^k + {\ov r}_{0,2}^k + {\ov r}_{1,2}^k$ for $k = 1,2$, so once
one solution $(\ov \th, {\ov \om}_0)$ to \tg{4-9} has been found, any
other can be written in the form $(\ov \th + \sum_k {\ov r}_{1,2}^k
\alpha_k,\ {\ov \om}_0 + \sum_k ({\ov r}_{2,0}^k + {\ov r}_{0,2}^k)
\alpha_k)$ for unique 1-forms $\alpha_k$, $k = 1,2$.

Now we define the 2-forms \[ \aligned \ov \Theta & = d\ov \th + \ppair
{{\ov \om}_0} {\ov \th} 1\\ \ov \Phi & = d{\ov \om}_0 - \frac12 \left(
\pair {{\ov \om}_{2,0}} {{\ov \om}_{2,0}} {1,0} + \pair {{\ov
\om}_{0,2}} {{\ov \om}_{0,2}} {0,1} \right) - \ov \Omega, \endaligned
\tag4-10$$ where $\ov \Omega$ is given by replacing $a_{i,j}$ and $\th$
by ${\ov a}_{i,j}$ and $\ov \th$ respectively in \tg{4-4}.

After some calculation, the exterior derivative of \tg{4-9} can be
written in the form \[ 0 = \ov J(\ov \Theta + \ov \Phi). \]

This implies that there are 2-forms ${\ov \Psi}_1$ and ${\ov \Psi}_2$
such that \[ \matrix \ov \Theta = \sum_k {\ov r}_{1,2}^k {\ov \Psi}_k,
& \text{and} & \ov \Phi = \sum_k ({\ov r}_{2,0}^k + {\ov r}_{0,2}^k)
{\ov \Psi}_k. \endmatrix \tag4-11$$

Substituting these relations into \tg{4-10} and differentiating, we
compute that \[ \matrix 0 = \sum_k {\ov r}_{1,2}^k d{\ov \Psi}_k, &
\text{and} & 0 = \sum_k ({\ov r}_{2,0}^k + {\ov r}_{0,2}^k) d{\ov
\Psi}_k. \endmatrix \]

Since the functions ${\ov r}_{1,2}^k + {\ov r}_{2,0}^k + {\ov
r}_{0,2}^k$ are linearly independent, we conclude that \[ \matrix d{\ov
\Psi}_k = 0 & \text{for $k = 1,2$.} \endmatrix \]

Let $U \subset \C$ be an open set on which the ${\ov \Psi}_k$'s are
{\it exact}. Clearly, $\C$ can be covered by such open sets. Let
$\alpha_k$ be 1-forms on $U$ such that $d\alpha_k = {\ov \Psi}_k$. If
we replace the pair $(\ov \th, \ov \om_0)$ by $\left(\ov \th - \sum_k
{\ov r}_{1,2}^k \alpha_k,\ {\ov \om}_0 - \sum_k \left({\ov r}_{2,0}^k +
{\ov r}_{0,2}^k\right) \alpha_k \right)$, then another calculation
shows that for this new pair \[ \aligned 0 =\ & d{\ov \th} + \ppair
{{\ov \om}_0} {{\ov \th}} 1\\ 0 =\ & d{\ov \om}_0 - \frac12 \left(
\pair {{\ov \om}_{2,0}} {{\ov \om}_{2,0}} {1,0} + \pair {{\ov
\om}_{0,2}} {{\ov \om}_{0,2}} {0,1} \right) - \ov \Omega. \endaligned
\tag4-12$$

Note that from \tg{4-10} and \tg{4-11} it follows that $\ov \Psi_k$ has
no $(0,2)$-part. Thus, $\alpha_k$ can be chosen to be of type $(1,0)$.
But then \tg{4-12} implies that $\ov \th$ and $\ov \om_0$ are
holomorphic 1-forms.

Now we let $F_0 := U \times \Bbb C^2$ with coordinates $(\ov a, \ov b,
s_1, s_2)$ and define the 1-forms \[ (\th, \om_0) := \left(\ov \th +
\sum_k {\ov r}_{1,2}^k ds_k,\ {\ov \om}_0 + \sum_k \left({\ov
r}_{2,0}^k + {\ov r}_{0,2}^k\right) ds_k \right). \]

Then it is not hard to show that $(\th, \om_0)$ is a holomorphic
coframe on $F_0$ satisfying the postulates of the Theorem.

The uniqueness of $(\th, \om_0)$ follows from the standard facts about
mappings preserving coframings \cite{G}. \qed
\enddemo

We are now ready to prove the existence result for
$H_{1,2}$-connections.

\proclaim{Corollary 4.8} For any constants $c, c_1, c_2 \in \Bbb C$ and
any point $u \in \C(c, c_1, c_2)$, there exists a regular torsion free
connection on some $H_{1,2}$-structure $\pi: F_0 \to \Z$ where $\Z$ is
some sixfold so that the image of the curvature map $K: F_0 \to \C_{c,
c_1, c_2}$ contains $u$.
\endproclaim

\demo{Proof} Let $U \subseteq \C_{c, c_1, c_2}$ be an open neighborhood
of $u$ for which the conclusion of Theorem 4.7. holds, i.e. there is a
principal $\Bbb C^2$-bundle $K: \tilde F_0 \to U$ and a coframe
$(\tilde \th, \tilde \om)$ on $\tilde F_0$ satisfying \tg{4-1} -
\tg{4-7} with $\om_{0,0} = 0$.

Pick a point $v \in \tilde F_0$ with $K(v) = u$. Since by the structure
equations we have $d\tilde \th \equiv 0 \mod \tilde \th$, it follows
that $\tilde F_0$ is foliated by integral leafs on which $\tilde \th$
vanishes. For some sufficiently small neighborhood $V$ of $v$, there
exists a submersion $\pi: V \to \Z$ onto some sixfold $\Z$ such that
$\ker(\pi_*) = \tilde \th^\perp$.

Moreover, standard arguments show that there is an inclusion $\imath: V
\hookrightarrow F_0 \subseteq \F$ of $V$ into an $H_{1,2}$-structure
$F_0$ on $\Z$ such that $\imath^*(\th) = \tilde \th$ where $\th$
denotes the tautological form on $F_0$.

Also, there is an unique $H_{1,2}$-connection $\om_0$ on $F_0$ with
$\imath^*(\om_0) = \tilde \om_0$. From the structure equations it is
then evident that the curvature map $K: F_0 \to \V$ satisfies $u \in
K(F_0) \subseteq \C_{c, c_1, c_2}$. \qed
\enddemo

As a consequence of the proof of Theorem 4.7. we have

\proclaim{Corollary 4.9} All regular torsion free $H_{1,2}$- and
$G_{1,2}$-connections are holomorphic. \qed
\endproclaim

\subheading{\S5 Summary}

In \S3 we have shown that the moduli space $\Z$ of rational curves $C$
in a complex threefold $\W$ whose normal bundle $N_C \to C$ is
equivalent to $\O(2) \oplus \O(2)$ forms a six dimensional manifold
which carries a natural $G_{1,2}$-structure $\pi: F \to \Z$. By Theorem
3.11. {\it most} of the intrinsic torsion of this structure vanishes. A
natural question is whether every holomorphic $G_{1,2}$-structure whose
torsion is of the form of Theorem 3.11. arises from such a moduli
space.

The answer is negative in general. The reason is that Proposition 2.9.
gives some first order restriction which is not automatically satisfied
if the torsion is of the form of Theorem 3.11., not even when the
torsion vanishes.

Before stating this result, let us write out the decompositions \[
\matrix \th & = & & \th_{1,2} & x \otimes x^2 & + & \th_{1,0} & x
\otimes x y & + & \th_{1,-2} x \otimes y^2\\ & & + & \th_{-1,2} & y
\otimes x^2 & + & \th_{-1,0} & y \otimes x y & + & \th_{-1,-2} y
\otimes y^2 \endmatrix \] and \[ \matrix \om & = & & \om_{0,0} & 1
\otimes 1\\ & & + & \om^{02}_{0,2} & 1 \otimes x^2 & + & \om^{02}_{0,0}
& 1 \otimes x y & + & \om^{02}_{0,-2} & 1 \otimes y^2\\ & & + &
\om^{20}_{2,0} & x^2 \otimes 1 & + & \om^{20}_{0,0} & x y \otimes 1 & +
& \om^{20}_{-2,0} & y^2 \otimes 1. \endmatrix \]

\proclaim{Proposition 5.1} Suppose $\Z$ is the moduli space of rational
curves in the threefold $\W$ whose normal bundle is equivalent to
$\O(2) \oplus \O(2)$, and suppose furthermore that the associated
$G_{1,2}$-structure $\pi: F \to \Z$ is torsion free. Then the torsion
free connection on $F$ is locally symmetric.
\endproclaim

\demo{Proof} By Theorem 3.11. the torsion free connection on $\pi: F
\to \Z$ must be special. Thus, by Proposition 2.9., \[ \Cal I :=
\ker(pr_1 \circ \pi_{F,\N})_* = \{ \th_{1,-2}, \th_{-1,-2},
\om^{02}_{0,-2} \}. \]

Of course, this means that $\Cal I$ must satisfy the {\it Frobenius
condition} $d\Cal I \equiv 0 \mod \Cal I$. However, from the structure
equations \tg{4-1} - \tg{4-4} we compute that $d\om^{02}_{0,-2} \equiv
9 \pair {a_{0,2}} {x^2} 2 \th_{1,0} \w \th_{-1,0} \mod \Cal I$. From
here it follows that the Frobenius condition is satisfied if and only
if $a_{0,2} \equiv 0$ on $F$. By \tg{4-5} and \tg{4-6}, this implies
that $\bb \equiv 0$ and $\pair {a_{2,0}} {a_{2,0}} {2,0} = \frac 34
c$.

However, $\bb$ represents the covariant derivative of the curvature
tensor. It follows that $F$ is locally symmetric. \qed
\enddemo

Let us now consider the question which regular torsion free
$G_{1,2}$-connections $(\pi: F \to \Z, \th, \om)$ admit a restriction
$(\pi: \hat F \to \Y, \hat \th, \hat \om)$ in the sense of
Definition~3.6. Since the holonomy of a torsion free
$G_{1,2}$-connection is contained in $H_{1,2}$ it follows that the
holonomy of the restriction to $\Y$ is contained in $H_3$.

\proclaim{Proposition 5.2} Let $(\th, \om)$ be a regular torsion free
$H_{1,2}$-connection on $\pi: F_0 \to \Z$ with structure constants $(c,
c_1, c_2)$. Then $F_0$ admits a restriction to an $H_3$-connection
$\pi: \hat F_0 \to \Y$ if and only if $c_1 = 0$.

In this case, the restriction $\hat F_0$ is uniquely determined, and
the connection on $\pi: \hat F_0 \to \Y$ is regular in the sense of
\cite{Br}.

Conversely, given a regular $H_3$-connection on $\pi: \hat F_0 \to \Y$,
there is a unique regular torsion free $H_{1,2}$-connection which
extends the connection on $\hat F_0$.
\endproclaim

\demo{Proof} Let $\pi: F_0 \to \Z$ be the torsion free regular
$H_{1,2}$-connection. If a restriction on $\hat F_0 \to \Y$ exists then
$T\hat F_0$ must be annihilated by the ideal \[ \Cal J = \{ \th_{-1,0}
- 2 \th_{1,-2},\ \th_{1,0} - 2 \th_{-1,2},\ \om^{02}_{0,i} -
\om^{20}_{i,0}, i = 0,1,2 \}. \]

Thus, $\Cal J$ must satisfy the {\it Frobenius condition} $d\Cal J
\equiv 0 \mod \Cal J$. A calculation using the structure equations
\tg{4-1} - \tg{4-4} yields that this is the case if and only if \[ 2
a_{2,0} = 3 a_{0,2}. \tag5-1$$

Taking the exterior derivative $2 da_{2,0} - 3 da_{0,2} \mod \Cal J$,
we conclude that $\bb$ must be of the form \[ \bb = x \otimes b^3_x + y
\otimes b^3_y \tag5-2$$ for some $\V_3$-valued function $b^3$ where the
subscripts stand for partial derivatives.

Let us define $\hat F_0 \subseteq F_0$ by \tg{5-1} and \tg{5-2}. Then
it is evident that any reduction of $F_0$ must be contained in $\hat
F_0$.

{}From the structure equations \tg{4-5} - \tg{4-7} we calculate that the
differentials of the components of \tg{5-1} and \tg{5-2} are linearly
independent. Also, substituting \tg{5-1} and \tg{5-2} into $f^c_1$ from
Proposition 4.2. we calculate $f^c_1 = 0$.

Therefore, $\hat F_0 = \emptyset$ if $c_1 \neq 0$. Conversely, if $c_1
= 0$ one can verify that $\hat F_0$ is non-empty and hence an eight
dimensional analytic submanifold of $F_0$. Moreover, $\dim(T\hat F_0
\cap \ker\pi_*) \equiv 4$, and so $\Y := \pi(\hat F_0)$ is an analytic
submainfold of $\Z$. Now it is easy to verify that $(\pi: \hat F_0 \to
\Y, \hat \th, \hat \om)$ with $\hat \th$ and $\hat \om$ as in
Definition 3.6. is the desired restriction. Of course, \tg{5-1} and
\tg{5-2} determine $\hat F_0$ uniquely.

Note that this restriction is a torsion free $H_3$-connection. The
final statement follows from the classification of regular
$H_3$-connections in \cite{Br}. They are uniquely determined by two
constant parameters, and it left to the reader to verify that these
correspond to the constants $c$ and $c_2$. \qed
\enddemo

\definition{Remark}
\roster

\item It seems likely that the last statement in Proposition 5.2. holds
true even if the $H_3$-connection on $\Y$ is {\it not} regular, i.e. in
this case there should still be an extension to a unique torsion free
$H_{1,2}$-connection. There does not seem to be any substantial
obstacle to proving this other than the immense calculations required
to determine the non-regular $H_{1,2}$-connections.

\item Since every torsion free holomorphic $G_3$\=connection is
equivalent to the moduli space of contact curves in a contact threefold
$\W$ \cite{Br}, it follows from the results in \S3 that {\it every
holomorphic torsion free $G_3$-connection can be extended to a
$G_{1,2}$-connection whose torsion is given as in the proof of} Theorem
3.13.

A characterization of $H_3$-connections is therefore that they are
precisely those $G_3$-connections which can be extended to a {\it
torsion free} $G_{1,2}$\=connection. (cf. Theorem 0.4.)

\item A somewhat surprising aspect comes from a comparison of
Proposition~5.1. and Proposition 5.2. Namely, if $\Y$ admits a regular
torsion free $H_3$\=connection then on the one hand, by \tg2, the
connection on $\Y$ can be extended to a connection on the moduli space
$\Z$ of rational curves in $\W$.

On the other hand, if we let $\Z'$ be the {\it torsion free} extension
from Proposition 5.2. then it follows from Proposition 5.1. that $\Z'$
is different from $\Z$ unless both $\Z$ and $\Y$ are flat: indeed, the
only locally symmetric $H_3$-connection is the flat one. (cf. Theorem
0.1.)

In other words, the extension $\Z$ of $\Y$ which seems most natural in
the geometric sense is different from the extension $\Z'$ of $\Y$ which
is most natural from the torsion point of view.

\endroster
\enddefinition

\definition{Definition 5.3} Let $\Cal P$ be a complex five dimensional
manifold. A {\it linear rank 2 Pfaffian system on $\Cal P$} or, for
short, a {\it Pfaffian structure} on $\Cal P$ is a differential ideal
$\Cal I$ on $\Cal P$ with the property that, locally, there is a
holomorphic coframe $\kappa_1, \kappa_2, \alpha, \beta_1, \beta_2$ on
$\Cal P$ such that \[ \Cal I = \{ \kappa_1, \kappa_2 \} \] and \[
d\kappa_i = \alpha \w \beta_i \mod \Cal I, \ \ i = 1,2. \]

A curve $C \subseteq \Cal P$ is called an {\it integral curve} if the
tangent vectors of $C$ are anihilated by $\Cal I$.
\enddefinition

A Pfaffian structure on $\Cal P$ may also be regarded as a rank 2
subbundle $L \subseteq T^*\Cal P$ where $L$ is locally spanned by
$\kappa_1$ and $\kappa_2$.

For example, if $\W$ is any three dimensional manifold then $\Cal P:=
\P{}T\W$ carries a canonical Pfaffian structure \cite{EDS}. Namely, for
local coordinates $(x,y,z)$ and $(x,y,z,u,v)$ on $\W$ and $\Cal P$
respectively such that the bundle map $\pi: \Cal P \to \W$ is given by
$(x,y,z,u,v) \mapsto (x,y,z)$, this system is given as $\Cal I := \{dy
- u dx, dz - v dx\}$.

Thus, for a curve of the form $(x, y(x), z(x))$ in $\W$ there is a
unique integral lift to $\Cal P$, namely $(x, y(x), z(x), y'(x),
z'(x))$.

A key observation is now given by the following

\proclaim{Theorem 5.4} Let $\pi: F \to \Z$ be a holomorphic
$G_{1,2}$-connection whose torsion is of the form of Theorem 3.11.
Then, locally, $\Z$ is (contained in) the moduli space of all integral
curves $C$ of a Pfaffian structure $\Cal I$ on some fivefold $\Cal P$.
\endproclaim

\demo{Proof} The proof requires to compute the structure equations for
connections whose torsion is of the required form. Since these
equations are quite complex and since we shall not need them any
further, they will be omitted.

It follows from these equations that the differential ideal \[ \Cal J
:= \{ \th_{1,0}, \th_{-1,0}, \th_{1,-2}, \th_{-1,-2}, \om^{02}_{0,-2}
\} \] {\it does} satisfy the Frobenius condition $d\Cal J \equiv 0 \mod
\Cal J$. Thus, at least locally, there is a map $p: F \to \Cal P$ onto
some five dimensional complex manifold $\Cal P$ such that $\ker(p_*) =
\Cal J^\perp$.

For each point $t \in \Z$, we let $C_t := p(\pi^{-1} (t))$. It is then
easy to see that $C_t$ is a rational curve in $\Cal P$, and hence we
may regard $\Z$ as the moduli space of certain rational curves in $\Cal
P$.

Let $\Cal I_0 := \{ \th_{1,-2}, \th_{-1,-2} \}$. Then a calculation
shows that for each vector field $X \in \Cal J^\perp$ on $F$, $\Cal L_X
(\Cal I_0) \subseteq \Cal I_0 \mod \Cal J$. Therefore, there is a
differential system $\Cal I$ on $\Cal P$ such that $p^*(\Cal I) = \Cal
I_0$.

Taking the exterior derivatives of $\th_{1,-2}$ and $\th_{-1,-2}$ it
follows that $\Cal I$ is indeed a Pfaffian structure on $\Cal P$.
Moreover, since $\pi^{-1} (t)$ is integral to $\Cal I_0$ for all $t \in
\Z$ it follows that $C_t$ is an integral curve for all $t \in \Z$. \qed
\enddemo

Theorem 5.4. suggests that it should be more natural to regard
$G_3$-structures, $H_3$\=structures and $G_{1,2}$-structures as moduli
spaces of integral curves of a fivefold with Pfaffian structure rather
than as curves in a threefold. Indeed, the remarks preceding Theorem
5.4. indicate how the moduli space of curves in a threefold $\W$ may be
regarded merely as a special case of this.

It should also be instructive to see how the local invariants of a
Pfaffian structure on $\Cal P$ \cite{C} relate to the associated
$G_{1,2}$-structure. This will be pursued in a sequel of the present
paper.

\enddocument